\documentclass[11pt]{amsart} 
 \usepackage{amsmath} 
 \usepackage{epsfig} 
 \usepackage{graphics} 
 
\def\bp{\mathbf p} 
\def\bq{\mathbf q} 
\def\bn{\mathbf n} 

\def\l{\langle} 
\def\r{\rangle} 
\def\Ga{\Gamma} 
\def\la{\lambda}

\def\p{\partial}

\newtheorem{Theorem}{Theorem}[section] 
\newtheorem{Corollary}{Corollary}[section] 
 
\newtheorem{Lemma}{Lemma}[section] 
\newtheorem{Proposition}{Proposition}[section]

\theoremstyle{definition} 
 
\newtheorem{Remark}{Remark}[section]

\catcode`@=11 \@addtoreset{equation}{chapter} 
 
\catcode`@=12 
 
\title[ ] 
{\bf \Large Symmetry and Asymmetry: the method of moving spheres} 
 
\author[]{} 
\thanks{} 
\begin{document} 
\maketitle 
 
\centerline{\scshape Qinian Jin 
\begin{footnote} 
{Department of Mathematics, The University of Texas at Austin, 
1 University  Station C1200, Austin,  TX 78712, 
qjin@math.utexas.edu.} 
\end{footnote} 
YanYan Li 
\begin{footnote} 
{Department of Mathematics, Rutgers University, 
110 Frelinghuysen Rd., 
Piscataway, NJ 08854, yyli@math.rutgers.edu. 
Partially supported by a NSF grant.} 
\end{footnote} 
and Haoyuan Xu 
\begin{footnote} 
{Department of Mathematics, Rutgers University, 
110 Frelinghuysen Rd., 
Piscataway, NJ 08854, hyxu@math.rutgers.edu.} 
\end{footnote} 
}

 

\baselineskip 14.3pt 
 
\section{\bf Introduction} 
\setcounter{equation}{0}

In this paper we will consider some nonlinear elliptic equations 
on ${\mathbb R}^n$ and ${\mathbb S}^n$. We first consider 
 
\begin{equation}\label{1.1} 
\Delta u+\frac{c}{|x|^2} u +u^{(n+2)/(n-2)}=0\quad 
\mbox{and}\quad  u>0 \mbox{ in } {\mathbb R}^n\backslash \{0\}. 
\end{equation} 
Through the work of Obata \cite{O71}, Gidas-Ni-Nirenberg 
\cite{GNN} and Caffarelli-Gidas-Spruck \cite{CGS}, the asymptotic 
behavior of solutions of (\ref{1.1}) as well as the classification of global solutions 
are well understood in the case when $c=0$. In \cite{V} V\'{e}ron 
raised the following question: {\it For $c\in 
{\mathbb R}$, $c\ne 0$ and $n\ge 3$, let $u\in C^\infty({\mathbb 
R}^n\backslash \{0\})$ satisfy (\ref{1.1}). Is it true that $u$ 
must be radially symmetric about the origin?} He pointed out 
that there might be non-radial solutions of certain form as suggested in section 4 of \cite{BV}. The 
following result partially answers this question. 
 
\begin{Theorem}\label{T1.1} 
\begin{enumerate} 
 
\item[(i)] If $c\ge (n-2)^2/4$, then (\ref{1.1}) has no 
smooth solution. 
 
\item[(ii)] If $c>0$ then any $u\in C^\infty({\mathbb R}^n\backslash 
\{0\})$ satisfying (\ref{1.1}) must be radially symmetric about 
the origin, $u'(r)<0$ for $0<r<\infty$, and there exists a 
positive constant $C$ such that $u(x)\le C|x|^{-(n-2)/2}$ 
for $x\in {\mathbb R}^n\backslash \{0\}$.

\item[(iii)] If $c<(n-2)^2/4$ then (\ref{1.1}) has 
infinitely many smooth radial solutions. 
 
\item[(iv)] For any $c<-(n-2)/4$, (\ref{1.1}) has non-radial 
solutions. Moreover, the number of non-radial solutions goes to $\infty$ when $c\to -\infty$. 
\end{enumerate} 
\end{Theorem} 
 
We remark that the non-radial solutions we produced for (iv) are not 
of the form suggested in \cite{BV} and it remains to be an interesting question to study the existence of solutions of the suggested form. We also remark that the above question of 
V\'{e}ron remains open for $-(n-2)/4\le c<0$. (iii) was a 
known result, a proof can be found in \cite{Fowler}. It is interesting to note that the number $(n-2)^2/4$ 
appeared in Theorem \ref{T1.1} is exactly the best constant in the 
classical Hardy inequality which states that 
$$ 
\int_{{\mathbb R}^n}|\nabla u|^2 dx\ge \left(\frac{n-2}{2}\right)^2 
\int_{{\mathbb R}^n} \frac{u^2}{|x|^2} dx 
$$ 
for all $u\in C_0^\infty({\mathbb R}^n)$ with $n\ge 3$. 
 
In order to prove the radial symmetry of a solution $u$ of (\ref{1.1}) with $c>0$, we will use the method of moving spheres, a variant of the method of moving planes \cite{GNN}, to compare $u$ with its Kelvin transforms $u_{\bar{x}, \la}$: 
$$ 
u_{\bar x, \la}(x):=\left(\frac{\la}{|x-\bar x|}\right)^{n-2}u\left(\bar x+\frac{\la^2(x-\bar x)}{|x-\bar x|^2}\right), \quad x\in {\mathbb R}^n\setminus\{\bar x\}, 
$$ 
where $\la>0$ and $\bar x\in {\mathbb R}^n$. 
 
In order to find non-radial solutions of (\ref{1.1}) for $c<-(n-2)/4$, let 
$v(t,\theta):=e^{-\frac{n-2}{2}t}u(r,\theta)$, where $(r,\theta)$, $0<r<\infty$, 
$\theta\in {\mathbb S}^{n-1}$,  are the polar coordinates of ${\mathbb R}^n$ and 
$t=-\log r$. Then $u$ is a solution of (\ref{1.1}) if and only if $v$ satisfies 
the equation 
$$ 
v_{tt}+\Delta_{{\mathbb S}^{n-1}}v 
+\left(c-\frac{(n-2)^2}{4}\right)v+v^{(n+2)/(n-2)}=0 \quad \mbox{on } {\mathbb S}^{n-1}, 
$$ 
where ${\mathbb S}^{n-1}$ is the unit sphere with the canonical 
metric $g_0$ induced from ${\mathbb R}^{n}$, and $\Delta_{{\mathbb 
S}^{n-1}}$ is the corresponding Laplace-Beltrami operator on ${\mathbb 
S}^{n-1}$. If $v$ depends only on $\theta\in {\mathbb S}^{n-1}$, then 
\begin{equation}\label{1.2} 
\Delta_{{\mathbb S}^{n-1}}v+\left(c-\frac{(n-2)^2}{4}\right)v 
+v^{(n+2)/(n-2)}=0 \quad \mbox{and} \quad v>0 \quad \mbox{on } 
{\mathbb S}^{n-1}. 
\end{equation} 
The way we prove Theorem \ref{T1.1} (iv) is to show the 
existence of non-constant solutions of (\ref{1.2}). 
 
Set $N=n-1$, we will consider the existence of non-constant 
solutions of the equation 
\begin{equation}\label{1.3} 
-\Delta_{{\mathbb S}^N} v=v^p- \lambda v \quad \mbox{and} \quad v>0 
\quad \mbox{on } {\mathbb S}^N, 
\end{equation} 
where $\lambda\in {\mathbb R}$, $N\ge 2$ and  $1<p<N^*$, here $N^*$ denotes 
$(N+2)/(N-2)$ if $N\ge 3$ and $\infty$ if $N=2$ respectively. It 
is clear that (\ref{1.3}) has no solution if $\lambda\le 0$. It 
was proved in \cite{BV}, which sharpened an earlier result in \cite{GS81}, that if $0<\lambda\le N/(p-1)$ 
then the only solution of (\ref{1.3}) is the constant 
$v=\lambda^{1/(p-1)}$. Since we will use this result, let us state 
it in the following form. 
 
\begin{Theorem} \label{T1.2} $($\cite{BV,GS81}$)$ 
If $1<p<N^*$, then the only solution of (\ref{1.3}) is the constant $v=\la^{1/(p-1)}$ for every $0<\la\le N/(p-1)$. 
\end{Theorem}

By using bifurcation theories and a priori estimates of 
solutions, we will show that (\ref{1.3}) has non-constant 
solutions for every $\lambda>N/(p-1)$. In \cite{BV} Bidaut-V\'{e}ron 
and V\'{e}ron showed that for $1<p<N^*$ and 
$\lambda>N/(p-1)$ but close to $N/(p-1)$ there 
exists non-constant solution of (\ref{1.3}) due to the local 
bifurcation theory. In \cite{BL06} Brezis and Li gave a somewhat 
different proof which also implies that for $p>N^*$ and 
$\lambda<N/(p-1)$ with $|\lambda-N/(p-1)|$ small, 
(\ref{1.3}) has non-constant solutions. See \cite{BP06} for some 
related works. Our result will employ the global 
bifurcation theorem of Rabinowitz \cite{R}. 
 
In order to state our result more precisely, let us 
introduce some notations and terminology. Let $O(N+1)$ be the 
group consisting of $(N+1)\times(N+1)$ orthogonal matrices, and 
let $G$ be the subgroup of $O(N+1)$ consisting of those elements 
which fix ${\mathbf e}_{N+1}=(0, \cdots, 0, 1)$. We say a function 
$v$ defined on ${\mathbb S}^N$ is $G$-invariant if 
$v(O\theta)=v(\theta)$ for $\theta\in {\mathbb S}^N$ and $O\in G$. 
It is clear that a $G$-invariant function $v$ on ${\mathbb S}^N$ 
can be written as $v(\theta)=\tilde{v}(\theta_{N+1})$ for some 
function $\tilde{v}$ defined on $[-1, 1]$, where $\theta_{N+1}$ 
denotes the $(N+1)$-th coordinate of $\theta$. In the following, 
given a $G$-invariant function $v$ on ${\mathbb S}^N$, we will 
always use $\tilde{v}$ to denote the corresponding function 
defined on $[-1, 1]$ such that 
$v(\theta)=\tilde{v}(\theta_{N+1})$. For each integer $l\ge 0$ and 
$0<\alpha<1$ we introduce the Banach space 
$$ 
C^{l, \alpha}_{G}({\mathbb S}^N)=\left\{v\in C^{l, \alpha}({\mathbb 
S}^N): v \mbox{ is } G \mbox{-invarant}\right\}, 
$$where $C^{l,\alpha}({\mathbb S}^N)$ denotes the usual H\"{o}lder spaces.

We will show that (\ref{1.3}) has non-constant solutions in $C^{2, 
\alpha}_{G}({\mathbb S}^N)$ when $\lambda>N/(p-1)$. For each $k$ we set 
\begin{align}\label{1.4} 
{\mathcal S_k}:=\big\{&v\in C^{2, \alpha}_{G}({\mathbb S}^N): 
\tilde{v} \mbox{ has exactly } k \mbox{ zeroes, all of}\nonumber\\ 
&\mbox{them are in } (-1, 1) \mbox{ and  are simple}\big\}. 
\end{align} 
Clearly ${S_j}$ are mutually disjoint. 
Our bifurcation result reads as follows. 
 
\begin{Theorem}\label{T1.3} 
Assume that $N\ge 2$ and $1<p<N^*$. Let 
$\nu_k=k(k+N-1)$ and $\lambda_k: =\nu_k/(p-1)$ for $k\ge 1$. Then for each $k\ge 1$ and 
$\lambda_k<\lambda\le\lambda_{k+1}$, (\ref{1.3}) has $k$ distinct non-constant solutions $v_1$, 
$\cdots$, $v_k$ such that $v_j-\la^{1/(p-1)}\in {\mathcal S}_j $ 
for $1\le j\le k$. 
\end{Theorem}

Next we will consider symmetry properties of solutions of some elliptic 
equations on ${\mathbb S}^n$. A point $\theta\in {\mathbb S}^n$ is represented as 
$\theta=(\theta_1,\cdots,\theta_{n+1})\in {\mathbb R}^{n+1}$ with $\sum \theta_i^2=1$. 
In the following we will always use $\bn$ and ${\mathbf s}$ to 
denote the north pole and south pole respectively, i.e. $\bn=(0, \cdots, 0, 1)$ and 
${\mathbf s}=(0, \cdots, 0, -1)$.  When $n\ge 3$, 
the conformal Laplacian on ${\mathbb S}^n$ is defined as 
$$ 
{\mathcal L}_{{\mathbb S}^n}:=\Delta_{{\mathbb S}^n}-\frac{n(n-2)}{4}. 
$$ 
In the following $g$ denotes a given function in $C^0(\Omega\times(0,\infty))$, where 
$\Omega$ is one of the sets ${\mathbb S}^n$, ${\mathbb S}^n\setminus\{\bf n\}$ 
or ${\mathbb S}^n\setminus\{\bf n,\bf s\}$. which should be clear from the context. 
 
We first consider the equation 
\begin{equation}\label{1.5} 
-{\mathcal L}_{{\mathbb S}^n}v=g(\theta, v) \quad \mbox{and} \quad v>0 \mbox{ 
on } {\mathbb S}^n\backslash\{{\mathbf n}\}. 
\end{equation} 
We will give the symmetry property of solutions of (\ref{1.5}) under various conditions on $g$. 
The following conditions are used in the first result. 
 
\begin{enumerate} 
\item[(g1)] For each $s>0$ the function $\theta\to g(\theta, s)$ is rotationally symmetric about 
the line through $\bn$ and ${\mathbf s}$, 
 
\item[(g2)] For each $s>0$ and any $\theta,\theta'$ on the same geodesic passing through $\bn$ and ${\mathbf s}$, the function $\theta\to g(\theta, s)$ satisfies $g(\theta, s)>g(\theta',s)$ if 
$\theta_{n+1} > \theta'_{n+1}$. 
 
\item[(g3)] For each $\theta\in{\mathbb S}^n\setminus\{\bn, {\bf s}\}$ the function $s\to s^{-(n+2)/(n-2)} g(\theta, s)$ is non-increasing on $(0, \infty)$. 
 
\item[(g4)] For each $\theta\in {\mathbb S}^n\setminus\{\bn, {\bf s}\}$ the function $s\to g(\theta, s)$ is non-decreasing on $(0, \infty)$. 
\end{enumerate}

\begin{Theorem}\label{T1.4} 
For $n\ge 3$, assume that $g$ is continuous on $({\mathbb S}^n\setminus\{\mathbf n\})\times(0,\infty)$ with $g(\cdot,s)$ bounded from below in ${\mathbb S}^n\setminus\{\mathbf n\}$ for each $s\in (0,\infty)$ and satisfies (g1)-(g4). If $v\in C^2({\mathbb S}^n\setminus\{\mathbf n\})$ is a solution of (\ref{1.5}) satisfying 
\begin{equation}\label{1.6} 
\liminf_{\theta\to{\bn}}v(\theta)>0, 
\end{equation} 
then $v$ is rotationally symmetric about the line through $\bn$ 
and ${\mathbf s}$. 
\end{Theorem} 
 
\begin{Remark}\label{R1.1} 
Condition (\ref{1.6}) ensures that our moving sphere procedure can 
start, the proof can be found in Lemma \ref{L3.1} in Section 3. If 
we assume that $g(\theta,v)\ge 0$ on ${\mathbb 
S}^n\setminus\{{\mathbf n}\}$, then (\ref{1.6}) is satisfied automatically (see 
Lemma \ref{L3.2}). 
\end{Remark} 
 
The proof of Theorem \ref{T1.4} is based on a moving sphere procedure on 
${\mathbb S}^n$, with a feature of varying both the radius and the center of the moving sphere, which we will introduce in the following. 
Given a function $v$ on ${\mathbb S}^n$, let us first define its Kelvin transforms. 
Fix $\bp\in {\mathbb S}^n$ and $0<\lambda<\pi$, let 
$B_\la(\bp)$ be the geodesic ball on ${\mathbb S}^n$ with center 
$\bp$ and radius $\lambda$. Set $\Sigma_{\bp, \lambda}:={\mathbb 
S}^n\setminus \overline{B_\lambda(\bp)}$. Let $\varphi_{\bp, \lambda}: 
{\mathbb S}^n\to {\mathbb S}^n$ be the uniquely determined conformal 
diffeomorphism such that $\varphi_{\bp, 
\lambda}(B_\la(\bp))=\Sigma_{\bp, \lambda}$, $\varphi(\Sigma_{\bp, 
\lambda})=B_\la(\bp)$ and $\varphi_{\bp, \la}$ fixes every point on $\p B_\la(\bp)$. 
Then the Kelvin transforms of $v$ are defined by 
\begin{equation}\label{1.7} 
v_{\bp, \la}:=|J_{\varphi_{\bp, \la}}|^{\frac{n-2}{2n}} (v\circ 
\varphi_{\bp, \la}), 
\end{equation} 
where $J_{\varphi_{\bp, \la}}$ denotes the Jacobian of 
$\varphi_{\bp, \la}$. By the conformal invariance we have 
\begin{equation}\label{1.8} 
-{\mathcal L}_{{\mathbb S}^n} v_{\bp, \la} =|J_{\varphi_{\bp, 
\la}}|^{\frac{n+2}{2n}} \left(-{\mathcal L}_{{\mathbb S}^n} 
v\right)\circ \varphi_{\bp, \la}. 
\end{equation} 
If we use $(r, \omega)$, $0<r<\pi$, $\omega\in {\mathbb S}^{n-1}$, to 
denote the geodesic polar coordinates on ${\mathbb S}^n$ with respect 
to $\bp$, then 
$$ 
\varphi_{\bp, \lambda}(r, \omega)=(h_\lambda(r), \omega), 
$$ 
where $h_\lambda(r)\in (0, \pi)$ is determined by the equation 
$$ 
\cos h_\lambda(r)=\frac{2 \cos \lambda-(1+\cos^2\lambda) \cos 
r}{1+\cos^2\lambda-2\cos \lambda \cos r}. 
$$ 
Some straightforward calculation then gives 
$$ 
|J_{\varphi_{\bp, \la}}|(r, 
\omega)=\left(\frac{\sin^2\la}{1+\cos^2\la-2\cos\la 
\cos r}\right)^n 
$$ 
>From this we can see that $|J_{\varphi_{\bp, \la}}|<1$ on 
$\Sigma_{\bp, \la}$ if $\la<\pi/2$ and $|J_{\varphi_{\bp, \la}}|>1$ on 
$\Sigma_{\bp, \la}$ if $\la>\pi/2$. 
 
For a solution $v$ of (\ref{1.5}), the proof of its rotational 
symmetry is reduced to showing that $v=v_{\bp,\pi/2}$ on $\Sigma_{\bp, \pi/2}\setminus \{{\mathbf n}\}$ for every 
$\bp\in \p B_{\bp,\pi/2}({\mathbf n})$. The comparison of $v$ with $v_{\bp, \la}$ is always possible for small $\la>0$ if $v$ is regular at $\bp$, i.e. there exists $\la_0>0$ such that 
$$ 
v_{\bp, \la}\le v\quad \mbox{on } \Sigma_{\bp, \la}\setminus\{{\mathbf n}\} \mbox{ for each } 0<\la<\la_0. 
$$ 
The number $\la_0$ in general depends on $\bp$. Under the conditions in Theorem \ref{T1.4} we can show that $\la_0$ can be taken as $\pi/2$ if $\bp={\mathbf s}$. 
We define 
$$ 
\Sigma:=\{\bp\in {\mathbb S}^n: v\ge v_{\bp, \pi/2} \mbox{ in } \Sigma_{\bp, \pi/2}\setminus\{{\mathbf n}\}\} 
$$ 
By using the strong maximum principle and the Hopf lemma we are able to show 
$\Sigma\supset \overline{B_{\pi/2}({\mathbf s})}$. This is enough for our purpose. The way we prove $\Sigma\supset\overline{B_{\pi/2}({\mathbf s})}$ is of some independent interest: 
For any point $\bp\in\partial B_{\pi/2}(\mathbf s)$, we construct $x\in C^1([0,1],{\mathbb S}^n)$, $\la\in C^0([0,1],[0,\pi/2])$, satisfying 
$$ 
x(0)={\mathbf s}, \,\,x(1)=\bp,\,\, \la(0)=0,\,\,\la(1)=\pi/2 
$$ 
and prove 
$$ 
v_{x(t),\la(t)}\le v \quad\mbox{on } \Sigma_{x(t),\la(t)}\setminus\{\mathbf n\} \mbox{ for all } 0\le t\le 1. 
$$ 
In fact, we take 
$$ 
x(t)={\mathbf s}\,\, \mbox{ and } \,\, \la(t)=t\pi \,\,\,\mbox{for }  0\le t\le \frac{1}{2} 
$$ 
and 
$$ 
\la(t)=\frac{\pi}{2}\quad\mbox{for}\quad \frac{1}{2}\le t\le 1, 
$$ 
while for $\frac{1}{2}\le t\le 1$, $x(t)$ goes from $\mathbf s$ to $\bp$ along the shortest geodesic (the largest circle).

We next give a symmetry result on the equation 
\begin{equation}\label{1.9} 
-{\mathcal L}_{{\mathbb S}^n} v=g(\theta, v) \quad \mbox{and} \quad v>0 \quad 
\mbox{on } {\mathbb S}^n\setminus\{{\mathbf n}, {\mathbf s}\}, 
\end{equation} 
where a solution $v$ is allowed to have two singularities. For the function $g$, 
in addition to (g1) and (g3), we will assume the following two conditions. 
 
\begin{enumerate} 
\item[(g5)] For each $s>0$ and any $\theta,\theta'$ on the same geodesic passing through $\bn$ and $\mathbf s$, 
there holds $g(\theta, s)\ge g(\theta', s)$ if $\theta_{n+1}> \theta'_{n+1}>0$ and $g(\theta, s)\ge g(\theta', s)$ if $\theta_{n+1}<\theta'_{n+1}<0$. 
 
\item[(g6)] Either the inequalities in (g5) are strict or the function in (g3) is strictly decreasing. 
\end{enumerate} 
 
\begin{Theorem}\label{T1.5} 
For $n\ge 3$, assume that $g$ is continuous on $({\mathbb S}^n\setminus\{\bn,{\mathbf s}\})\times(0,\infty)$ with $g(\cdot,s)$ bounded from below in ${\mathbb S}^n\setminus\{\mathbf n,\mathbf s\}$ for each $s\in (0,\infty)$ and satisfies (g1), (g3), (g5) and (g6). If 
$v$ is a solution of (\ref{1.9}) on ${\mathbb S}^n\setminus\{\bn,{\mathbf s}\}$ satisfying 
\begin{equation}\label{1.10} 
\liminf_{\theta\to{\bn}}v(\theta)>0\quad \mbox{and} \quad \liminf_{\theta\to{\mathbf s}}v(\theta)>0, 
\end{equation} 
then $v$ is rotationally symmetric about the line through $\bn$ and ${\mathbf s}$. 
\end{Theorem} 
 
\begin{Remark}\label{R1.2} 
Similar to Remark \ref{R1.1}, (\ref{1.10}) is sufficient for the moving sphere procedure to start on the equator. 
If we assume that $g(\theta,v)\ge 
0$ on ${\mathbb S}^n\setminus\{{\mathbf n}, {\mathbf s}\}$, then (\ref{1.10}) is automatically satisfied. 
\end{Remark} 
 
As the first application of Theorem \ref{T1.4} and Theorem \ref{T1.5}, we consider the Matukuma equation 
\begin{equation}\label{1.11} 
-\Delta u=\frac{1}{1+|x|^2}u^p \quad \mbox{and} \quad  u>0\quad\mbox{in } {\mathbb R}^n, 
\end{equation} 
where $n\ge 3$ and $p\ge 0$. Let $\pi_{\mathbf n}: {\mathbb S}^n\setminus\{{\mathbf n}\}\to {\mathbb R}^n$ be the stereographic projection which sends ${\mathbf n}$ to $\infty$. 
Let $g_0$ be the standard metric on ${\mathbb S}^n$. It is well-known that 
$$ 
(\pi_{\mathbf n}^{-1})^*(g_0)=\xi(x)^{4/(n-2)} \sum_{i=1}^n (dx_i)^2, 
$$ 
where 
$$ 
\xi(x)=\left(\frac{2}{1+|x|^2}\right)^{(n-2)/2}, \quad x\in {\mathbb R}^n. 
$$ 
For a solution $u$ of (\ref{1.11}), we define a function $v$ on ${\mathbb S}^n\setminus \{{\mathbf n}\}$ by 
$$ 
v(\theta)=(\xi^{-1} u)(\pi_{\mathbf n}(\theta)), \quad \theta\in {\mathbb S}^n\setminus\{{\mathbf n}\}. 
$$ 
By using the conformal invariance one can check $v$ satisfies (\ref{1.5}) with 
$$ 
g(\theta, s):=\frac{1}{2} (1-\theta_{n+1})^{\frac{n-2}{2}(p-\frac{n}{n-2})} s^p. 
$$ 
Thus $g$ satisfies (g1)--(g4) if $0\le p<\frac{n}{n-2}$ and $g$ satisfies (g1), (g3), (g5) and (g6) if $p=\frac{n}{n-2}$. Theorem \ref{T1.4} and Theorem \ref{T1.5} then imply that $v$ is rotationally symmetric about the line through $\bn$ 
and ${\mathbf s}$, which in turn implies that $u$ is radially symmetric about the origin. On the other hand, it is easy to see that there is no smooth positive radially symmetric solutions to $-\Delta u\ge \frac{1}{1+|x|^2}u^p$ in ${\mathbb R}^n$ for $0\le p<1$. 
We thus obtain 
 
\begin{Corollary}\label{C1.1} 
If $1\le p\le n/(n-2)$, then any smooth solution of the Matukuma equation (\ref{1.11}) must be radially symmetric about the origin. If $0\le p<1$, (\ref{1.11}) has no smooth solutions. 
\end{Corollary} 
 
\begin{Remark}\label{R1.3} 
The result, as far as we know, is new for the case $n\ge 4$ and for the case $n=3$ and $0\le p\le 1$. When $n=3$ and $1<p<5$,  
the result was proved by Li in \cite{LY}. For $n\ge 3$ and 
$1<p<\frac{n+2}{n-2}$, the result was proved earlier by Li and Ni 
in \cite{LN1,LN2,LN3} under an additional finite total mass 
condition: $\int_{{\mathbb R}^n}\frac{1}{1+|x|^2}u^pdx<\infty$. 
Our method is different from theirs. In \cite{LN1,LN2,LN3} under 
the finite total mass condition, they analyzed the asymptotic 
behaviors of solutions at $\infty$ to ensure that the moving plane 
method can start at $\infty$. In \cite{LY}, Li obtained the 
asymptotic behavior of the solution in dimension $n=3$ without the 
finite total mass condition. Our proof, via the method of moving spheres, does not need to analyze the asymptotic behavior of 
solutions at $\infty$. 
\end{Remark} 
 
We now consider a special form of equation (\ref{1.9}) as follows 
\begin{equation}\label{1.12} 
-{\mathcal L}_{{\mathbb S}^n} v=K(\theta) v^{(n+2)/(n-2)} \quad \mbox{and} \quad v>0 \quad 
\mbox{on } {\mathbb S}^n\setminus\{{\mathbf n}, {\mathbf s}\}, 
\end{equation} 
i.e. $g(\theta,s)=K(\theta)s^{(n+2)/(n-2)}$, where $K(\theta)$ is a function defined on ${\mathbb S}^n \setminus\{\bn, {\mathbf s}\}$. As a consequence of Theorem \ref{T1.5} we have 
 
\begin{Corollary}\label{C1.2} 
For $n\ge 3$, assume that $K$ is continuous and 
on ${\mathbb S}^n\setminus\{{\mathbf n},{\mathbf s}\}$  and that $K$ is 
rotationally symmetric about the line through $\bn$ and ${\mathbf s}$. 
Assume further that there exists $-1<c<1$ such that $K(\theta)>K(\theta')$ if $\theta_{n+1}>\theta'_{n+1}\ge c$ 
and $K(\theta)<K(\theta')$ if $\theta_{n+1}<\theta'_{n+1}\le c$. Then any solution 
$v\in C^2({\mathbb S}^n\setminus\{{\mathbf n}, {\mathbf s}\})$ of (\ref{1.12}) satisfying 
(\ref{1.10}) is 
rotationally symmetric about the line through $\bn$ and ${\mathbf s}$. 
\end{Corollary} 
 
To see this, let $\varphi_c: {\mathbb S}^n\to {\mathbb S}^n$ be the conformal diffeomorphism 
such that $\varphi_c(B_{\pi/2}({\bf s}))=\Gamma_c$ and $\varphi_c(\Sigma_{\bf s, \pi/2})= 
{\mathbb S}^n\setminus \Gamma_c$, where $\Gamma_c:=\{\theta\in {\mathbb S}^n: \theta_{n+1}<c\}$. 
For a solution $v$ of (\ref{1.12}), we define 
$$ 
\hat{v}:=|J_{\varphi_c}|^{2n/(n-2)} (v\circ \varphi_c), 
$$ 
where $J_{\varphi_c}$ is the Jacobian of $\varphi_c$. Then 
$$ 
-{\mathcal L}_{{\mathbb S}^n} \hat{v}=(K\circ \varphi_c) \hat{v}^{(n+2)/(n-2)} \quad \mbox{and} 
\quad \hat{v}>0 \quad \mbox{on } {\mathbb S}^n\setminus \{\bn, {\mathbf s}\}. 
$$ 
We can use Theorem \ref{T1.5} to conclude that $\hat{v}$ is rotationally symmetric about the line 
through $\bn$ and ${\mathbf s}$ and so is $v$.

Theorem \ref{T1.5} can be used to classify $C^2$ solutions of the 
equation 
\begin{equation}\label{1.13} 
-{\mathcal L}_{{\mathbb S}^n} v=f(v) \quad \mbox{and} \quad v>0 \quad 
\mbox{on } {\mathbb S}^n, 
\end{equation} 
where $f: [0, \infty)\to {\mathbb R}$ is a continuous function. 
 
\begin{Corollary}\label{C1.3} 
Assume that $n\ge 3$ and that $s^{-(n+2)/(n-2)}f(s)$ is strictly decreasing 
on $(0, \infty)$. Then any solution $v\in 
C^2({\mathbb S}^n)$ of (\ref{1.13}) must be a constant $v\equiv c$ on 
${\mathbb S}^n$ satisfying $f(c)=\frac{n(n-2)}{4} c$. 
\end{Corollary} 
 
When $f$ satisfies some differentiability condition, Gidas and Spruck obtained this result in \cite{GS81}  by using Obata type argument. 
On the other hand Brezis and Li obtained it in \cite{BL06} by transforming the equation into an equation in ${\mathbb R}^n$ and then using the result established in \cite{GNN1} by moving plane method.

If we take 
\begin{equation}\label{1.14} 
g(\theta, s)=s^{\frac{n+2}{n-2}} +\left(\frac{n(n-2)}{4}-\beta\right)s 
\end{equation} 
for some number $\beta\in {\mathbb R}^n$, then (\ref{1.9}) reduces to 
the form 
\begin{equation}\label{1.15} 
-\Delta_{{\mathbb S}^n} v+\beta v=v^{\frac{n+2}{n-2}}\quad \mbox{and} 
\quad v>0 \quad \mbox{on } {\mathbb S}^n\backslash \{\bn, {\mathbf s}\}. 
\end{equation} 
By using Theorem \ref{T1.5} we can analyze the solutions of (\ref{1.15}) in some detail. 
 
\begin{Corollary}\label{C1.4} 
\begin{enumerate} 
\item[(i)] If $\beta<n(n-2)/4$ then any solution $v\in 
C^2({\mathbb S}^n\setminus\{\bn, {\mathbf s}\})$ of (\ref{1.15}) is 
rotationally symmetric about the line through $\bn$ and ${\mathbf s}$. 
 
\item[(ii)] If $\beta\le 0$ then (\ref{1.15}) has no solution 
$v\in C^2({\mathbb S}^n\setminus\{\bn, {\mathbf s}\})$. 
 
\item[(iii)] If $(n-2)/2<\beta<n(n-2)/4$, then the 
equation (\ref{1.15}) has infinitely many solutions in 
$C^2({\mathbb S}^n\setminus\{\bn, {\mathbf s}\})$, which have exactly tow singularities. 
\end{enumerate} 
\end{Corollary} 
 
It is well known that if $v\in C^\infty({\mathbb 
S}^n\setminus \{\bn, {\mathbf s}\})$ is a solution of (\ref{1.15}) satisfying 
$v\in L^{\frac{2n}{n-2}}({\mathbb S}^n)$, then $v\in C^\infty({\mathbb S}^n)$ and 
thus $v$ must be constant if $0<\beta<n(n-2)/4$.  Corollary \ref{C1.4} (iii) indicates that 
(\ref{1.15}) has non-constant solutions at least for $(n-2)/2<\beta<n(n-2)/4$ if one drops the condition 
$v\in L^{\frac{2n}{n-2}}({\mathbb S}^n)$. 
It is interesting to point out that, for $\beta=\beta_0:=\frac{1}{16}(n-2)(3n-2)$, 
$$ 
v(\theta):=\left(\frac{n-2}{2}\right)^{\frac{n-2}{2}} 
(1-\theta_{n+1})^{-\frac{n-2}{4}} 
$$ 
is also a solution of (\ref{1.15}) which has only one singularity at $\bn$. Note that $\beta_0\le (n-2)/2$ 
for $n=3, 4$, therefore it is probably true that (\ref{1.15}) has 
many solutions even if $0<\beta\le (n-2)/2$.

Similar problems can be considered on ${\mathbb S}^2$. We first 
consider the equation of the form 
\begin{equation}\label{1.16} 
-\Delta_{{\mathbb S}^2}v+1 =K(\theta)e^{2v}+f(\theta) 
\quad\mbox{on } {\mathbb S}^2\setminus\{\bn,{\mathbf s}\}, 
\end{equation} 
where $K$ and $f$ are continuous functions on ${\mathbb S}^2$ 
satisfying the following conditions. 
 
\begin{enumerate} 
\item[(K1)] For any $\theta,\theta'$ on the same geodesic passing through $\bn$ and ${\mathbf s}$, 
$K(\theta)\ge K(\theta')$ if $\theta_{n+1}>\theta'_{n+1}\ge 0$ and $K(\theta)\le 
K(\theta')$ if $\theta_{n+1}<\theta'_{n+1}\le 0$. 
 
\item[(f1)] For any $\theta,\theta'$ on the same geodesic passing through $\bn$ and ${\mathbf s}$, 
$f(\theta)\ge f(\theta')$ if $\theta_{n+1}>\theta'_{n+1}\ge 0$ and $f(\theta)\le 
f(\theta')$ if $\theta_{n+1}<\theta'_{n+1}\le 0$. 
 
\item[(Kf1)] $f^2+(K--\hspace{-0.35cm}\int_{{\mathbb S}^2}K)^2$ is not identically 
zero on ${\mathbb S}^2$. 
\end{enumerate}

Similar to Theorem \ref{T1.5} we have 
 
\begin{Theorem}\label{T1.6} 
Assume that $K$ and $f$ are continuous non-negative functions 
defined on ${\mathbb S}^2\setminus\{\bn, {\mathbf s}\}$ satisfying 
(K1), (f1) and (Kf1). If both $K$ and $f$ are 
rotationally symmetric about the line through $\bn$ and ${\mathbf 
s}$, then any solution $v\in C^2({\mathbb 
S}^2\setminus\{\bn,\mathbf s\})$ of (\ref{1.16}) satisfying 
\begin{equation}\label{1.17} 
\limsup_{\theta\rightarrow \bn}\frac{v(\theta)}{\log 
d(\theta,\bn)}\le 0 \quad \mbox{and} \quad 
\quad\limsup_{\theta\rightarrow {\mathbf s}}\frac{v(\theta)}{\log 
d(\theta, {\mathbf s})}\le 0 
\end{equation} 
must be rotationally symmetric about the line through $\bn$ and 
${\mathbf s}$. 
\end{Theorem} 
 
\begin{Remark}\label{R1.4} 
It is well known that if $f^2+(K--\hspace{-0.35cm}\int_{{\mathbb S}^2}K)^2\equiv 0$ on ${\mathbb S}^2$, the conclusion of Theorem \ref{T1.6} is not true. 
\end{Remark} 
 
As an immediate consequence of Theorem \ref{T1.6} we have 
 
\begin{Corollary}\label{C1.5} 
Suppose that $v\in C^2({\mathbb S}^2)$ satisfies the equation 
\begin{equation}\label{1.18} 
-\Delta_{{\mathbb S}^2}v+\beta =e^{2v} \quad\mbox{on } {\mathbb 
S}^2 
\end{equation} 
with $0<\beta<1$. Then $v$ must be constant. 
\end{Corollary}

Next we consider the equation 
\begin{equation}\label{1.19} 
-\Delta_{{\mathbb S}^2}v+1 =K(\theta)e^{2v}+f(\theta) 
\quad\mbox{on } {\mathbb S}^2\setminus\{\bn\}, 
\end{equation} 
where $K$ and $f$ are non-negative continuous functions on ${\mathbb S}^2$ 
satisfying the following conditions. 
 
\begin{enumerate} 
\item[(K2)] For any $\theta,\theta'$ on the same geodesic passing through $\bn$ and ${\mathbf s}$, 
$K(\theta)\ge K(\theta')$ if $\theta_{n+1}>\theta'_{n+1}$. 
 
\item[(f2)]  For any $\theta,\theta'$ on the same geodesic passing through $\bn$ and 
${\mathbf s}$, $f(\theta)\ge f(\theta')$ if $\theta_{n+1}>\theta'_{n+1}$. 
 
\item[(Kf2)]  $(f--\hspace{-0.35cm}\int_{{\mathbb S}^2}f)^2+(K--\hspace{-0.35cm}\int_{{\mathbb 
S}^2}K)^2$ is not identically zero on ${\mathbb S}^2$. 
\end{enumerate} 
 
Similar to Theorem \ref{T1.4} we have 
 
\begin{Theorem}\label{T1.7} 
Assume that $K$ and $f$ are continuous non-negative functions 
defined on ${\mathbb S}^2\setminus\{\bn\}$ satisfying (K2),(f2) 
and (Kf2). If both $K$ and $f$ are rotationally 
symmetric about the line through $\bn$ and ${\mathbf s}$, then any 
solution $v\in C^2({\mathbb S}^2\setminus\{\bn\})$ of (\ref{1.19}) 
satisfying 
\begin{equation}\label{1.20} 
\limsup_{\theta\rightarrow \bn}\frac{v(\theta)}{\log 
d(\theta,\bn)}\le 0 
\end{equation} must be rotationally symmetric about 
the line through $\bn$ and ${\mathbf s}$. 
\end{Theorem} 
 
We give an application of Theorem \ref{T1.7} to the mean 
field equation 
\begin{equation}\label{1.21} 
\Delta_{{\mathbb S}^2} \varphi +\frac{\exp(\alpha 
\varphi(\theta)-\gamma\l{\mathbf s}, \theta\r)}{\int_{{\mathbb 
S}^2} \exp(\alpha\varphi(\sigma)-\gamma\l {\mathbf s}, \sigma\r) 
d\sigma}-\frac{1}{4\pi}=0 \quad \mbox{on } {\mathbb S}^2, 
\end{equation} 
where $\alpha$ and $\gamma$ are nonnegative numbers, It is easy to 
see that the new function 
$$ 
v=\frac{\alpha}{2}\varphi(\theta)-\frac{1}{2}\log \int_{{\mathbb 
S}^2} \exp(\alpha \varphi(\sigma)-\gamma\l {\mathbf s}, \sigma\r) 
d\sigma +\frac{1}{2} \log \frac{\alpha}{2} 
$$ 
satisfies the equation 
$$ 
-\Delta_{{\mathbb S}^2} v+\frac{\alpha}{8\pi} =\exp(2v-\gamma\l 
{\mathbf s}, \theta\r) \quad \mbox{on } {\mathbb S}^2. 
$$ 
This is exactly the equation (\ref{1.20}) with 
$K(\theta)=\exp(-\gamma\l {\mathbf s}, \theta\r)$ and 
$f(\theta)=1-\frac{\alpha}{8\pi}$. By using Theorem \ref{T1.7} we 
thus conclude the following result which was proved by Lin 
\cite{LCS}. 
 
\begin{Corollary}\label{C1.6} 
If $0<\alpha<8 \pi$ and $\gamma\ge 0$, then any solution of 
(\ref{1.21}) must be rotationally symmetric about the line 
through $\bn$ and ${\mathbf s}$. 
\end{Corollary}

Using the moving sphere procedure, we can also show the following result. 
We consider the equation 
\begin{equation}\label{1.22} 
-{\mathcal L}_{{\mathbb S}^n} v=K(\theta) v^{(n+2)/(n-2)} \quad \mbox{and}\quad v>0 \quad 
\mbox{on } {\mathbb S}^n\setminus\{\bn\}. 
\end{equation} 
We will use $\nabla_X$ to denote the covariant differentiation on ${\mathbb 
S}^n$ with respect to the vector field $X$. 
 
\begin{Theorem}\label{T1.8} 
Assume that $n\ge 3$ and that $K$ is a $C^1$ 
function on ${\mathbb S}^n$ such that 
$\nabla_{\frac{\partial}{\partial \theta_{n+1}}}K\ge 0$ and 
is not identically zero on ${\mathbb S}^n\setminus\{\bn, {\bf s}\}$. 
If $v\in C^2({\mathbb S}^n\setminus\{\bn\})$ is a 
solution of (\ref{1.22}) satisfying 
\begin{equation}\label{1.23} 
\liminf_{\theta\to{\mathbf n}} v(\theta)>0, 
\end{equation} 
then 
\begin{equation}\label{1.24} 
v(\theta)\ge C_0 d(\theta, \bn)^{-(n-2)/2} \quad \mbox{for all } \theta\in {\mathbb S}^n\setminus\{\bn\} 
\end{equation} 
for some positive constant $C_0$, where $d(\cdot, \cdot)$ denotes 
the distance on ${\mathbb S}^n$. In particular, (\ref{1.22}) has no 
$C^2({\mathbb S}^n)$ solutions. 
\end{Theorem} 
 
\begin{Remark}\label{R1.5} 
If $K\ge 0$, then (\ref{1.23}) is automatically satisfied. 
\end{Remark} 
 
\begin{Remark}\label{R1.6} 
The non-existence of $C^2({\mathbb S}^n)$ solution to (\ref{1.22}) is due to Kazdan-Warner \cite{KW}. 
\end{Remark}

Our moving sphere procedure can also be used to obtain a Kazdan-Warner type obstruction for some 
fully nonlinear elliptic equations on ${\mathbb S}^n$ for $n\ge 3$. 
 
 Let $g_0$ be the standard metric on ${\mathbb S}^n$. 
For $n\ge 3$, let $A_g$ denote the Schouten tensor of a metric $g$ 
\begin{equation}\label{1.25} 
A_g=\frac{1}{n-2}\left(Ric_g-\frac{R_g}{2(n-1)}g\right) 
\end{equation} 
where $Ric_g$ and $R_g$ denote the Ricci tensor and the scalar 
curvature of $g$ respectively. 
 
For $0<v\in C^2({\mathbb S}^n)$, we consider the conformal change of 
metric $g_1=v^{\frac{4}{n-2}}g_0$. Then 
\begin{align}\label{1.26} 
A_{g_1}=&-\frac{2}{n-2}v^{-1}\nabla_{g_0}^2 
v+\frac{2n}{(n-2)^2}v^{-2}\nabla_{g_0} 
v\otimes\nabla_{g_0}v\nonumber\\ 
&-\frac{2}{(n-2)^2}v^{-2}|\nabla_{g_0} v|^2 
g_0+A_{g_0}. 
\end{align}

We assume that 
\begin{equation}\label{1.27} 
 \Gamma\subset {\mathbb R}^n \mbox{ is an open convex symmetric cone with vertex at the origin} 
\end{equation} 
such that 
\begin{equation}\label{1.28} 
\Gamma_n\subset\Gamma\subset \Gamma_1 
\end{equation} 
with $\Gamma_1:=\{\lambda\in {\mathbb R}^n:\sum_{i}\lambda_i>0\}$ and 
$\Gamma_n:=\{\lambda\in {\mathbb R}^n:\lambda_i>0 \mbox{ for all } i\}$, 
where $\Gamma$ being symmetric means 
$(\lambda_1,...,\lambda_n)\in \Gamma$ implies 
$(\lambda_{i_1},...,\lambda_{i_n})\in\Gamma$ for any permutation 
$(i_1,...,i_n)$ of $(1,2,...,n)$. 
 
We also assume that $f$ is a function defined on $\Gamma$ such that 
\begin{equation}\label{1.29} 
f\in C^1({\Gamma})\mbox{ is symmetric in } \Gamma 
\end{equation} 
and 
\begin{equation}\label{1.30} 
  f>0 \mbox{ and } f_{\lambda_i}>0 \mbox{ in } \Gamma. 
\end{equation} 
 
Given a positive $C^1$ function $K$ on ${\mathbb S}^n$, we consider the equation 
\begin{equation}\label{1.31} 
f(\lambda(A_{g_1}))=K,\quad \lambda(A_{g_1})\in\Gamma \quad\mbox{and}\quad v>0\quad\mbox{on}\quad {\mathbb S}^n. 
\end{equation} 
where $g_1=v^{\frac{4}{n-2}}g_0$ and $\la(A_{g_1})$ denotes the eigenvalues of $A_{g_1}$ with respect to $g_1$.

\begin{Theorem}\label{T1.9} 
For $n\ge 3$, assume that $(f,\Gamma)$ satisfies condition 
(\ref{1.27})-(\ref{1.30}), and that $K$ 
is a positive $C^1$ function on ${\mathbb S}^n$, such that 
$\nabla_{\frac{\partial}{\partial \theta_{n+1}}}K\ge 0$ 
and is not identically zero on ${\mathbb S}^n\setminus\{\bn, {\bf s}\}$. 
Then (\ref{1.31}) has no $C^2({\mathbb S}^n)$ solutions. 
\end{Theorem}

\begin{Remark}\label{R1.7} 
The non-existence of $C^2({\mathbb S}^n)$ solution of (\ref{1.31})
under the assumption on $K$, or some similar ones,  was known for 
$(f,\Gamma)=(\sigma_k^{\frac{1}{k}},\Gamma_k)$, where 
$$ 
\sigma_k(\la)=\sum_{1\le i_1<\cdots<i_k\le n}\la_{i_1}\cdots\la_{i_k} 
$$ and 
$$ 
\Gamma_k=\{\la\in{\mathbb R}^n:\sigma_1(\la)>0,\cdots,\sigma_k(\la)>0\}. 
$$ 
For details, see the work of Viaclovsky \cite{Via00c}, Han \cite{Han06} and Delano\"{e} \cite{Del06}. 
Our method of proof is completely different from theirs. 
\end{Remark} 
 
We can also consider the following equation on the half sphere 
${\mathbb S}^n_+:=\{\theta\in{\mathbb 
S}^n: \theta_1\ge 0\}$ for $n\ge 3$: 
\begin{equation}\label{1.32} 
\left\{ 
\begin{array}{l} 
f\left(\lambda(A_{v^{\frac{4}{n-2}}g_0})\right)=K(\theta),\,\,\, \la(A_{v^{\frac{4}{n-2}}g_0})\in\Gamma 
\,\,\, \mbox{and}\, \,\,  v>0 \,\,\mbox{ on } {\mathbb S}^n_+ \\ 
\\ 
\frac{\partial 
v}{\partial\nu}=H(\theta)v^{\frac{n}{n-2}}\quad \mbox{on } 
\partial {\mathbb S}^n_+, 
\end{array}\right. 
\end{equation} 
where $\nu$ denotes the unit outer normal of 
$\partial {\mathbb S}^n_+$. 
 
\begin{Theorem}\label{T1.10} 
For $n\ge 3$, assume that $(f,\Gamma)$ satisfies 
(\ref{1.27})-(\ref{1.30}). Assume also that $K>0$, $K\in C^1({\mathbb S}^n_+)$ and 
$H\in C^1(\partial {\mathbb S}^n_+)$ such that $\nabla_{\frac{\partial}{\partial \theta_{n+1}}}K\ge 0$ on 
${\mathbb S}^n_+$, $\nabla_{\frac{\partial}{\partial \theta_{n+1}}}H\ge 0$ on $\p {\mathbb S}^n_+$ and at least one 
of these two inequalities is strict somewhere. Then (\ref{1.32}) has no positive $C^2({\mathbb S}^n_+)$ solutions. 
\end{Theorem} 
 
\begin{Remark}\label{R1.8} 
In the case that $f(\la)=\la_1+\cdots+\la_n$, the non-existence of $C^2({\mathbb S}^n_+)$ solution 
was known, see Bianchi and Pan \cite{GP}. Our proof, similar to that of Theorems \ref{T1.8}--\ref{T1.9}, 
is completely different. 
\end{Remark}

The paper is organized as follows. In section 2, we prove Theorem \ref{T1.1} and Theorem \ref{T1.3}. In section 3, we use 
moving sphere procedure on sphere to show the other theorems and 
their corollaries.

\section{\bf Some results on V\'{e}ron's question} 
\setcounter{equation}{0} 
 
In this section we give the proof of Theorem \ref{T1.1}. The 
radial symmetry of solutions of (\ref{1.1}) will be proved in 
subsection \ref{sub2.1} by using the method of moving spheres. The 
existence of non-radial solutions  of (\ref{1.1}) follows from 
Theorem \ref{T1.3} whose proof is based on a global bifurcation 
analysis and will be provided in subsection \ref{sub2.2}.

\subsection{Proof of Theorem \ref{T1.1}}\label{sub2.1} 
 
We first state a calculus lemma due to \cite{L}, which gives the 
symmetric property of a function through the investigation of its 
Kelvin transforms. 
 
\begin{Lemma}\label{L2.1} 
If $u\in C^1({\mathbb R}^n\backslash\{0\})$ is a function such that 
for each $y\ne 0$ there holds 
\begin{equation}\label{2.1} 
u_{y, \la}(x)\le u(x) \quad \forall 0<\la<|y| \mbox{ and } 
|x-y|\ge \la \mbox{ with } x\ne 0, 
\end{equation} 
then $u$ must be radially symmetric about the origin, and $u'(r)\le 0$ for 
$0<r<\infty$. 
\end{Lemma} 
\begin{proof} 
We include here the proof for completeness.  For any $x\in {\mathbb 
R}^n\backslash \{0\}$ and any number $a>0$, let ${\mathbf e}$ be any unit 
vector in ${\mathbb R}^n$ such that $\l x-a{\mathbf e}, {\mathbf e}\r <0$. For any number 
$\tau>a$, if we set $\la=\tau-a$ and $y=\tau {\mathbf e}$, then $0<\la <|y|$ 
and $|x-y|>\la$. So we may apply (\ref{2.1}) to get 
\begin{equation}\label{2.2} 
u(x)\ge \left(\frac{\tau-a}{|x-\tau {\mathbf e}|}\right)^{n-2} u\left(\tau {\mathbf e} 
+\frac{(\tau-a)^2(x-\tau {\mathbf e})}{|x-\tau {\mathbf e}|^2}\right). 
\end{equation} 
It is easy to check that 
$$ 
\tau {\mathbf e} +\frac{(\tau-a)^2(x-\tau {\mathbf e})}{|x-\tau {\mathbf e}|^2}\rightarrow 
x-2\left(\l x, {\mathbf e}\r -a\right){\mathbf e} \quad \mbox{as } \tau\rightarrow 
\infty. 
$$ 
Since $\l x-a {\mathbf e}, {\mathbf e}\r<0$, we must have $x-2\left(\l x, {\mathbf e}\r 
-a\right){\mathbf e}\ne 0$. Therefore, by sending $\tau\rightarrow \infty$ 
in (\ref{2.2}) and using the continuity of $u$ in ${\mathbb 
R}^n\backslash \{0\}$, we obtain 
$$ 
u(x)\ge u\left(x-2\left(\l x, {\mathbf e}\r -a\right){\mathbf e}\right). 
$$ 
This immediately implies that $u$ is radially symmetric about the 
origin and $u'(r)\le 0$ for $0<r<\infty$. 
\end{proof}

Instead of proving the symmetric property about solutions of 
(\ref{1.1}) directly, we consider the following  more general 
equation 
\begin{equation}\label{2.3} 
\Delta u + a(x) u+u^{(n+2)/(n-2)}=0 \quad \mbox{and} \quad u>0 
\mbox{ in } {\mathbb R}^n\backslash \{0\}, 
\end{equation} 
where $n\ge 3$, $\Delta$ denotes the Laplace operator on ${\mathbb 
R}^n$,  and $a: {\mathbb R}^n\backslash \{0\}\to [0, \infty)$ is a 
continuous function verifying the following property: 
 
\begin{enumerate} 
\item[(A)] for each $x\ne 0$ there holds 
$$ 
\left(\frac{\la}{|z|}\right)^4 a \left(x+\frac{\la^2z}{|z|^2} 
\right) < a(x+z), \quad \forall 0<\la <|x| \mbox{ and } 
|z|>\la. 
$$ 
\end{enumerate}

We have the following symmetry result for solutions of 
(\ref{2.3}).

\begin{Proposition}\label{P2.1} 
Let $a: {\mathbb R}^n\backslash \{0\}\to [0, \infty)$ be a 
continuous function satisfying (A). If $u\in C^2({\mathbb 
R}^n\backslash \{0\})$ is a solution of (\ref{2.3}), then $u$ must 
be radially symmetric about the origin and $u'(r)<0$ for all 
$0<r<\infty$. 
\end{Proposition} 
 
\begin{proof} The proof is based on 
the method of moving spheres. From (\ref{2.3}) it follows that 
$\Delta u\le 0$ and $u>0$ in ${\mathbb R}^n\backslash \{0\}$. So, by 
the maximum principle, we have 
\begin{equation}\label{2.4} 
\liminf_{|x|\rightarrow 0} u(x)>0 \quad \mbox{and} \quad 
\liminf_{|x|\rightarrow \infty} |x|^{n-2} u(x)>0. 
\end{equation} 
One can follow the proof of \cite[Lemma 2.1]{LZ} to conclude that 
for each $y\ne 0$ there exists $\la(y)>0$ such that 
$$ 
u_{y,\la}(x)\le u(x), \quad \forall 0<\la<\la(y) \mbox{ and } 
|x-y|\ge \la \mbox{ with } x\ne 0. 
$$ 
Define 
$$ 
\bar{\la}(y):=\left\{0<\mu\le |y|: u_{y, \la}(x)\le u(x), 
\forall 0<\la<\mu,\, |x-y|>\la \mbox{ with } x\ne 0\right\}. 
$$ 
Then $\bar{\la}(y)>0$. By using Lemma \ref{L2.1}, it suffices to 
show that 
\begin{equation}\label{2.5} 
\bar{\la}(y)= |y|, \quad \forall y\in {\mathbb 
R}^n\backslash\{0\}. 
\end{equation} 
Suppose that (\ref{2.5}) is not true, then there exists $y_0\ne 0$ 
such that $\bar{\la}(y_0)<|y_0|$. Let $\la_0:=\bar{\la}(y_0)$, 
then from the definition of $\la_0$ we have 
\begin{equation}\label{2.6} 
u_{y_0, \la_0}(x)\le u(x) \qquad \forall |x-y_0|>\la_0 \mbox{ with 
} x\ne 0. 
\end{equation} 
A straightforward calculation shows that 
$$ 
\Delta u_{y_0, \la_0}(x)=-\left(\frac{\la_0}{|x-y_0|}\right)^4 
\varphi\left(y_0+\frac{\la_0^2(x-y_0)}{|x-y_0|^2}\right) u_{y_0, 
\la_0}(x)- u_{y_0, \la_0}(x)^{(n+2)/(n-2)}. 
$$ 
Therefore, by using (A) and (\ref{2.6}), we have for $|x-y_0|> 
\la_0$ with $x\ne 0$ that 
\begin{align*} 
-\Delta u_{y_0, \la_0}(x)<&\varphi(x) u_{y_0, \la_0}(x) +u_{y_0, 
\la_0}(x)^{(n+2)/(n-2)}\le-\Delta u(x). 
\end{align*} 
This, together with the strong maximum principle and the Hopf 
lemma, gives 
\begin{equation}\label{2.7} 
u(x)-u_{y_0, \la_0}(x)>0, \quad \forall |x-y_0|>\la_0 \mbox{ with 
} x\ne 0, 
\end{equation} 
\begin{equation}\label{2.8} 
\liminf_{|x|\rightarrow 0} (u(x)-u_{y_0,\la_0}(x))>0 \,\,\mbox{ 
and }\,\, \liminf_{|x|\rightarrow \infty} |x|^{2-n}(u(x)-u_{y_0, 
\la_0}(x))>0, 
\end{equation} 
and 
\begin{equation}\label{2.9} 
\left.\frac{d}{dr}(u-u_{y_0, \la_0})\right|_{\p B_{\la_0}(y_0)}>0. 
\end{equation} 
Properties (\ref{2.7})-(\ref{2.9}) lead to, as in section 2 of 
\cite{LZ}, a contradiction to the definition of $\bar{\la}(y_0)$. 
For reader's convenience, we include a proof. From (\ref{2.9}) it 
follows that there exists $R_0$ satisfying $\la_0<R_0<|y_0|$ such 
that 
$$ 
\frac{d}{dr}(u-u_{y_0, \la})(x)>0 \quad \mbox{for } \la_0\le \la\le 
R_0 \mbox{ and } \la\le |x-y_0|\le R_0. 
$$ 
Since $u-u_{y_0, \la}=0$ on $\p B_\la(y_0)$, we have 
\begin{equation}\label{2.10} 
u(x)-u_{y_0, \la}(x)>0 \quad \mbox{ for } \la_0\le \la< R_0 \mbox{ 
and } \la<|x-y_0|\le R_0. 
\end{equation} 
>From (\ref{2.8}) one can find $c>0$, $R_1>|y_0|>R_0$ and $\eta>0$ 
such that 
$$ 
u(x)-u_{y_0,\la_0}(x)\ge\left\{ \begin{array}{lll} c|x-y_0|^{2-n} 
&& \mbox{ if } |x-y_0|\ge R_1\\ 
\\ 
c && \mbox{ if } 0<|x|<\eta. 
\end{array}\right. 
$$ 
But it easy to see that there exists $\varepsilon_1>0$ such that 
if $\la_0\le \la\le \la_0+\varepsilon_1$ then 
$$ 
|u_{y_0, \la}(x)-u_{y_0, \la_0}(x)|\le \left\{ 
\begin{array}{lll} 
\frac{1}{2}c|x-y_0|^{2-n} && \mbox{ if } |x-y_0|\ge R_1\\ 
\\ 
\frac{1}{2}c && \mbox{ if } 0< |x|\le \eta. 
\end{array}\right. 
$$ 
Therefore for $\la_0\le \la\le \la_0+\varepsilon_1$ there holds 
\begin{equation}\label{2.11} 
u(x)-u_{y_0, \la}(x)>0 \quad \mbox{ if } |x-y_0|\ge R_1 \mbox{ or 
} 0<|x|\le \eta. 
\end{equation} 
Finally, by continuity, (\ref{2.7}) implies that there exists 
$\varepsilon_2>0$ such that 
\begin{equation}\label{2.12} 
u(x)-u_{y_0, \la}(x)>0 \quad \mbox{if } \la_0\le \la\le 
\la_0+\varepsilon_2, \, R_0\le |x-y_0|\le R_1 \mbox{ and } |x|\ge 
\eta. 
\end{equation} 
Combining (\ref{2.10}), (\ref{2.11}) and (\ref{2.12}) we have for 
some $\varepsilon>0$ that 
$$ 
u(x)-u_{y_0, \la}(x)>0 \quad \mbox{if } \la_0\le \la\le 
\la_0+\varepsilon \mbox{ and } |x-y_0|\ge \la \mbox{ with } x\ne 
0. 
$$ 
This gives a contradiction to the definition of $\la_0$. We thus 
obtain (\ref{2.5}). 
\end{proof}

Now we are ready to give the proof of Theorem \ref{T1.1} by assuming 
Theorem \ref{T1.3}. 
 
\begin{proof}[Proof of Theorem \ref{T1.1}] 
We first show part (ii) by using Proposition \ref{P2.1}. Let 
$\varphi(x)=c/|x|^2$, then it suffices to verify (A) for 
$\varphi$. This is equivalent to showing that 
$$ 
\la^4|x+z|^2<|z|^4 \left|x+\frac{\la^2 z}{|z|^2}\right|^2, \quad 
\forall 0<\la<|x| \mbox{ and } |z|>\la. 
$$ 
It can be confirmed by the following computation 
\begin{align*} 
\la^4|x+z|^2-|z|^4 \left|x+\frac{\la^2 z}{|z|^2}\right|^2 
&=\left(\la^2-|z|^2\right)\left\{(\la^2+|z|^2)|x|^2+2\la^2\l x, 
z\r\right\}\\ 
&\le(\la^2-|z|^2)\left\{(\la^2+|z|^2)|x|^2-\la^2(|x|^2+|z|^2)\right\}\\ 
&= (\la^2-|z|^2)(|x|^2-\la^2)|z|^2 < 0. 
\end{align*} 
Thus Proposition \ref{P2.1} applies to conclude that $u$ must be 
radially symmetric about the origin, and consequently $u$ satisfies the 
ordinary differential equation 
$$ 
u''+\frac{n-1}{r} u'+\frac{c}{r^2} u +u^{(n+2)/(n-2)}=0,\quad 
u>0 \quad \mbox{and} \quad u'< 0 \,\, \mbox{for } r\in(0, \infty). 
$$ 
Then 
$$ 
(r^{n-1}u')'< -r^{n-1}u^{(n+2)/(n-2)}. 
$$ 
Therefore for any $r>\varepsilon>0$ we have 
$$ 
r^{n-1} u'(r)-\varepsilon^{n-1} u'(\varepsilon) 
<-\int_\varepsilon^r s^{n-1} u(s)^{(n+2)/(n-2)} ds\le 
-\frac{1}{n}(r^n-\varepsilon^n) u(r)^{(n+2)/(n-2)}. 
$$ 
Since $u'(\varepsilon)\le 0$ we may drop the second term on the 
left hand and then take $\varepsilon\rightarrow 0$ to get 
$$ 
r^{n-1}u'(r)\le -\frac{1}{n}r^n u(r)^{(n+2)/(n-2)} \quad 
\mbox{for all } r>0. 
$$ 
This is equivalent to 
$\frac{d}{dr}\left(u(r)^{-4/(n-2)}\right)\ge 
\frac{4}{n(n-2)} r$. Hence for $r\ge\varepsilon>0$ there holds 
$$ 
u(r)^{-4/(n-2)}\ge 
u(\varepsilon)^{-4/(n-2)}+\frac{2}{n(n-2)} 
\left(r^2-\varepsilon^2\right)\ge 
\frac{2}{n(n-2)}\left(r^2-\varepsilon^2\right). 
$$ 
Letting $\varepsilon\rightarrow 0$ we then obtain 
$$ 
u(r)\le C r^{-(n-2)/2} \quad \mbox{ for all } r> 0. 
$$ 
This gives the desired estimate. 
 
We now use (ii) to show (i). Suppose (\ref{1.11}) has a solution $u$ for some 
$c\ge (n-2)^2/4$. We define 
$$ 
w(t):=e^{-\frac{n-2}{2}t} u(e^{-t}), \quad t\in (-\infty, \infty). 
$$ 
One can verify that $w$ satisfies the ordinary differential 
equation 
\begin{equation}\label{2.13} 
w''+\left(c-\frac{(n-2)^2}{4}\right)w+w^{(n+2)/(n-2)}=0 \quad \mbox{and}\quad 
w>0 \mbox{  on } (-\infty, \infty). 
\end{equation} 
Thus $w$ is a positive strictly concave function defined on $(-\infty, \infty)$. However, such function does not exist. 
 
(iv) From Theorem \ref{T1.3} it follows that (\ref{1.2}) has 
non-constant solutions when $c<-(n-2)/4$, which give non-radial 
solutions of (\ref{1.1}). 
\end{proof}

\subsection{Proof of Theorem \ref{T1.3}: A bifurcation analysis}\label{sub2.2}

We first give a fact concerning solutions of (\ref{1.3}). 
 
\begin{Lemma}\label{L2.2} 
Suppose $1<p<N^*$. Then for any $\Lambda>0$ there exists a 
positive constant $C(N, p, \Lambda)$ depending only on $N$, $p$ 
and $\Lambda$ such that any non-constant solution $v$ of 
(\ref{1.3}) with $\lambda\le \Lambda$ satisfies 
$$ 
1/C(N, p, \Lambda)\le v\le C(N, p, \Lambda) \quad \mbox{on } {\mathbb 
S}^n. 
$$ 
\end{Lemma} 
 
\begin{proof} 
Note that (\ref{1.3}) has no solution for $\lambda\le 0$. Considering Theorem 
\ref{T1.2}, we may assume $N/(p-1)\le \la<\Lambda$. The upper bound 
can be obtained by using 
the blow-up technique together with the fact that the equation 
$-\Delta u=u^p$ in ${\mathbb R}^N$ has no positive solution if 
$1<p<N^*$ (see e.g. \cite{SY}). In order to get the lower 
bound, we first use the Harnack inequality to get 
$$ 
\max_{{\mathbb S}^N} v\le C_1(N, p, \Lambda) \min_{{\mathbb S}^n} v 
$$ 
for some positive constant $C_1(N, p,  \Lambda)$ depending only on 
$N$, $p$ and $\Lambda$. So it suffices to derive a lower bound of 
$\max_{{\mathbb S}^n} v$. Suppose $\max_{{\mathbb S}^n} v=v(x_0)$ for 
some $x_0\in {\mathbb S}^N$. Then $-\Delta_{{\mathbb S}^N} v(x_0)\ge 0$ 
and hence $v(x)^p-\lambda v(x_0)\ge 0$. This implies $v(x_0)\ge 
\lambda^{1/(p-1)}\ge (N/(p-1))^{1/(p-1)}$. 
\end{proof}

\begin{Lemma}\label{L2.3} 
The eigenvalues of $-\Delta_{{\mathbb S}^N}$ restricted to $C^{4, 
\alpha}_{G}({\mathbb S}^N)$ are $\nu_k=k(N+k-1)$ for $k=0, 
1,\cdots$, they are all simple, the eigenspace of $\nu_k$ is 
spanned by a function $p_k$ which can be written as 
$p_k(\theta)=\tilde{p}_k(\theta_{n+1})$, where $\tilde{p}_k(t)$ is 
a polynomial of degree $k$. Moreover, all the zeroes of 
$\tilde{p}_k(t)$ are simple and in $(-1, 1)$. 
\end{Lemma} 
 
\begin{proof} All the assertions can be found in \cite{BG} except 
the last part. In the following we will show by induction on $k$ 
that $\tilde{p}_k$ has exactly $k$ simple zeroes in $(-1, 1)$. 
This is clear for $k=0$ since $\tilde{p}_0=1$. Now we assume that 
$\tilde{p}_k$ has $k$ simple zeroes in $(-1, 1)$, say 
$-1<t_1<t_2<...<t_k<1$. Set $t_0=-1$ and $t_{k+1}=1$. It suffices 
to show that $\tilde{p}_{k+1}$ has a zero in each interval $(t_i, 
t_{i+1})$ for $i=0,\cdots,k$. Suppose for some $0\le i\le k$ the 
polynomial $\tilde{p}_{k+1}$ has no zeroes in $(t_i, t_{i+1})$, 
then both $\tilde{p}_k$ and $\tilde{p}_{k+1}$ do not change sign 
in this interval. Without loss of generality, we may assume both 
$\tilde{p}_k$ and $\tilde{p}_{k+1}$ are positive in $(t_i, 
t_{i+1})$. Let 
$$ 
\Sigma_i:=\left\{\begin{array}{lll} \{\theta\in {\mathbb S}^N: 
\theta_{N+1}<t_1\}, & \mbox{if } i=0,\\ 
\\ 
\{\theta\in {\mathbb S}^N: t_i<\theta_{N+1}<t_{i+1}\}, & \mbox{if } 
1\le i\le k-1,\\ 
\\ 
\{\theta\in {\mathbb S}^N: \theta_{N+1}>t_k\}, & \mbox{if } i=k. 
\end{array}\right. 
$$ 
Then both $p_k$ and $p_{k+1}$ are non-negative on the domain 
$\Sigma_i$, $p_k=0$ and $\frac{\p p_{k}}{\p \nu}\le 0$ on $\p 
\Sigma_i$. Recall that 
$$ 
-\Delta_{{\mathbb S}^N} p_k=\nu_k p_k \quad \mbox{and} \quad 
-\Delta_{{\mathbb S}^N} p_{k+1}=\nu_{k+1} p_{k+1} \quad \mbox{on } 
{\mathbb S}^n. 
$$ 
We have 
\begin{align*} 
(\nu_{k+1}-\nu_k)\int_{\Sigma_i} p_k p_{k+1} =& -\int_{\Sigma_i} 
p_k \Delta_{{\mathbb S}^N} p_{k+1} 
+\int_{\Sigma_i} p_{k+1} \Delta_{{\mathbb S}^N} p_k \\ 
=& \int_{\p \Sigma_i} \left(p_{k+1}\frac{\p p_k}{\p \nu}- p_k 
\frac{\p p_{k+1}}{\p \nu}\right)\le 0, 
\end{align*} 
which is a contraction. 
\end{proof} 
 
In order to study the solutions of (\ref{1.3}), let 
$v=\lambda^{1/(p-1)}(w+1)$, then $w$ satisfies the equation 
\begin{equation}\label{2.14} 
-\Delta_{{\mathbb S}^N} w=\lambda\left((w+1)^p-w-1\right) \quad 
\mbox{and} \quad w>-1 \quad \mbox{on } {\mathbb S}^N. 
\end{equation} 
We are going to find $G$-invariant non-zero solutions of 
(\ref{2.14}) for each $\lambda>N/(p-1)$.

\begin{Lemma}\label{L2.4} 
(a) For any $0<\lambda \le N/(p-1)$ the only solution of 
(\ref{2.14}) is $w=0$. 
 
(b) For any $\Lambda>0$, there exist positive constant $C(N, p, 
\Lambda)$ and $\varepsilon(N, p, \Lambda)$ depending only on $N$, 
$p$ and $\Lambda$ such that 
$$ 
-1+\varepsilon(N, p, \Lambda)\le w\le C(N, p, \Lambda) \quad 
\mbox{on } {\mathbb S}^N 
$$ 
for any solution $w$ of (\ref{2.14}) with $\lambda\le \Lambda$. 
 
(c) Any non-zero $G$-invariant solution $w$ of (\ref{2.14}) neither vanishes 
at the north pole nor at the south pole on ${\mathbb S}^n$; moreover, by 
writing $w(\theta)=\tilde{w}(\theta_{n+1})$ for some function 
$\tilde{w}$ on $[-1,1]$, then all zeroes of $\tilde{w}$ are 
in $(-1, 1)$ and are simple. 
\end{Lemma} 
 
\begin{proof} 
(a) follows from Theorem \ref{T1.2} and (b) is the consequence of 
Lemma \ref{L2.2}. In the following we will prove (c). Let us first 
show that $w$ does not vanish at the north pole. If it vanishes at 
the north pole, then the strong maximum principle implies that the 
north pole must be an accumulating point of zeroes of $w$ on 
${\mathbb S}^N$. Since $w$ is $G$-invariant, this would imply that 
all derivatives of $w$ at the north pole are zero. The unique 
continuation property then implies $w=0$ on ${\mathbb S}^N$. The 
same argument gives $w\ne 0$ at the south pole. Therefore all 
zeroes of $\tilde{w}$ are in $(-1, 1)$. Use (\ref{2.14}) one can 
see that $\tilde{w}$ satisfies the ordinary differential equation 
$$ 
-(1-t^2) \tilde{w}''+nt \tilde{w}'=\lambda 
\left((\tilde{w}+1)^p-\tilde{w}-1\right) \quad \mbox{on } [-1,1]. 
$$ 
Since $\tilde{w}\ne 0$ on $[-1,1]$, this implies that all zeroes 
of $\tilde{w}$ must be simple. 
\end{proof}

\begin{proof}[Proof of Theorem \ref{T1.3}] 
We will use the bifurcation theory to carry out the proof. We 
first formulate (\ref{2.14}) as an operator equation. Since 
$\Delta_{{\mathbb S}^N}$ is $O(N+1)$-invariant, it follows from the 
theory of elliptic equations that $-\Delta_{{\mathbb S}^N} +I: C^{4, 
\alpha}_{G}({\mathbb S}^N)\to C^{2, \alpha}_{G}({\mathbb 
S}^N)$ is invertible. Let $T$ denote its inverse, then $T: C^{2, 
\alpha}_{G}({\mathbb S}^N)\to C^{2, \alpha}_{G}({\mathbb 
S}^N)$ is a compact linear operator. Let 
$$ 
{\mathcal D}:=\left\{(w, \mu)\in C^{2, \alpha}_{G}({\mathbb 
S}^N)\times {\mathbb R}: \mu>1 \mbox{ and } w>-1 \mbox{ on } {\mathbb 
S}^N\right\} 
$$ 
and set 
$$ 
\mu=(p-1)\lambda+1\quad \mbox{and} \quad g(w, 
\mu):=\frac{\mu-1}{p-1} T\left((w+1)^p-pw-1 \right). 
$$ 
Then, finding a $G$-invariant nonzero solution of (\ref{2.14}) 
is equivalent to finding a nonzero solution of the operator 
equation 
\begin{equation}\label{2.15} 
 f(w, \mu):=w-\mu T w-g(w, \mu)=0 \quad \mbox{in } 
{\mathcal D}. 
\end{equation} 
It is clear that $g$ is a nonlinear compact map of ${\mathcal D}$ 
into $C^{2, \alpha}_{G}({\mathbb S}^N)$, and $g(w, 
\mu)=o(\|w\|_{C^{2,\alpha}({\mathbb S}^N)})$ uniformly on bounded 
$\mu$ interval. 
 
Note that $1/\mu$ is an eigenvalue of $T$ if and only if $\mu-1$ 
is an eigenvalue of $-\Delta_{{\mathbb S}^n}$ restricted to $C^{4, 
\alpha}_{G}({\mathbb S}^N)$. Therefore from Lemma \ref{L2.3} we 
can see that all the eigenvalues of $T$ are simple and given by 
$1/\mu_k$, where $\mu_k=1+k(N+k-1)$. It then follows from 
Krasnoselski's Theorem (see \cite[Theorem 3.3.1]{N}) that each 
$(0, \mu_k)$ is a bifurcation point of $f(w, \mu)=0$ in ${\mathcal 
D}$. 
 
Let $\mathcal{S}$ denote the closure of nontrivial solutions $(w, 
\mu)$ of $f(w, \mu)=0$ in ${\mathcal D}$. According to Lemma 
\ref{L2.4} (a), ${\mathcal S}\cap (C^{2, \alpha}_{G}({\mathbb 
S}^N)\times (1, \mu_1))=\emptyset$. Let ${\mathcal C}_k$ be the 
connected component of ${\mathcal S}$ containing $(0, \mu_k)$. 
Then by the global bifurcation theory of Rabinowitz (see 
\cite[Theorem 1.4]{R} or \cite[Theorem 3.4.1]{N}) we know that 
either (i) ${\mathcal C}_k$ is not compact in ${\mathcal D}$ or 
(ii) ${\mathcal C}_k$ contains a point $(0, \mu_j)$ with $j\ne k$. 
We are going to rule out the case (ii). 
 
For each $k\ge 1$ we define ${\mathcal S}_k$ as in 
(\ref{1.4}), i.e. 
\begin{align*} 
{\mathcal S_k}:=\big\{&v\in C^{2, \alpha}_{G}({\mathbb S}^N): 
\tilde{v} \mbox{ has exactly } k \mbox{ zeroes, all of}\nonumber\\ 
&\mbox{them are in } (-1, 1) \mbox{ and are simple} \big\}. 
\end{align*} 
It is clear that each ${\mathcal S}_k$ is an open set in $C^{2, 
\alpha}_G({\mathbb S}^N)$. By the local bifurcation theorem of 
Crandall and Rabinowitz (see \cite[Theorem 1.7]{CR}), near each 
bifurcation point $(0, \mu_k)$, ${\mathcal S}$ has the 
parametrization $(w_k(s), \mu_k(s))$, $|s|<a_k$ for some small 
$a_k>0$, where $\mu_k(0)=\mu_k$, $w_k(s)=s p_k+ s\psi_k(s)$ and 
$\psi_k(0)=0$. According to Lemma \ref{L2.3}, $p_k\in {\mathcal 
S}_k$, thus $w_k(s)\in {\mathcal S}_k$ for small $s\ne 0$. 
Therefore there exists a neighborhood ${\mathcal O}_k$ of $(0, 
\mu_k)$ in ${\mathcal D}$ such that if $(w, \mu)\in {\mathcal 
O}_k\cap {\mathcal S}$ and $w\ne 0$, then $w\in {\mathcal S}_k$. 
Let 
$$ 
{\mathcal B}_k:=\left\{(w, \mu)\in {\mathcal C}_k: w\in {\mathcal 
S}_k\right\}\cup \{(0, \mu_k)\}. 
$$ 
If we can show that $\mathcal{C}_k={\mathcal B}_k$ for each $k\ge 
1$, then ${\mathcal C}_k$ can not contain a point $(0, \mu_j)$ with $j\ne k$, 
and we therefore rule out the case (ii). 
 
In order to show ${\mathcal C}_k={\mathcal B}_k$, it suffices to 
show that ${\mathcal B}_k$ is both open and closed in ${\mathcal 
C}_k$. It is clear that ${\mathcal B}_k$ is open in ${\mathcal 
C}_k$. Suppose now that $\{(w^{(l)}, \mu^{(l)}\}$ is a sequence in 
${\mathcal B}_k$ such that $(w^{(l)}, \mu^{(l)})\rightarrow (w, 
\mu)$ in ${\mathcal C}_k$. Note that $(w, \mu)$ is a solution 
of (\ref{2.15}). If $w=0$ on ${\mathbb S}^N$, then $\mu=\mu_j$ for 
some $j$. If $j=k$, $(0,\mu_k)\in {\mathcal B}_k$; if $j\ne k$, then $w^{(l)}\in {\mathcal S}_j\cap 
{\mathcal S}_k$ for large $l$ which is impossible. Thus $w$ is a 
nonzero solution of (\ref{2.14}), (c) of Lemma \ref{L2.4} then 
implies $w\in {\mathcal S}_i$ for some $i$. If $i\ne k$, then the 
openness of ${\mathcal S}_i$ implies $w^{(l)}\in {\mathcal 
S}_i\cap {\mathcal S}_k$ for large $l$ which is again impossible. 
Hence $w\in {\mathcal S}_k$ and $(w, \mu)\in {\mathcal B}_k$. 
Therefore ${\mathcal B}_k$ is closed in ${\mathcal C}_k$. 
 
The above argument has ruled out case (ii), therefore each ${\mathcal 
C}_k$ is noncompact in ${\mathcal D}$. 
Let 
$$ 
\Lambda=\sup \left\{\lambda>\mu_k: {\mathcal C}_k\cap \left(C^{2, \alpha}_{G}({\mathbb S}^N)\times 
\{\mu\}\right)\ne \emptyset,\forall \mu_k\le \mu<\lambda\right\} 
$$ 
We will show that $\Lambda=\infty$. If not, say $\Lambda<\infty$. By connectedness of ${\mathcal C}_k$ and (a) of Lemma \ref{L2.4}, ${\mathcal C}_k\subset \left(C^{2, \alpha}_{G}({\mathbb S}^N)\times 
[\mu_1,\Lambda]\right)$. It follows, using (b) of Lemma \ref{L2.4} that ${\mathcal C}_k$ is compact, 
a contradiction. So $\Lambda=\infty$, i.e. 
 
$$ 
{\mathcal C}_k\cap \left(C^{2, \alpha}_{G}({\mathbb S}^N)\times 
\{\mu\}\right)\ne \emptyset 
$$ 
for any $\mu>\mu_k$. The proof is thus complete. 
\end{proof}

\section{\bf Symmetric results on ${\mathbb S}^n$} 
\setcounter{equation}{0} 
 
\subsection{Some preliminary results} 
 
Given a function $v$ on ${\mathbb S}^n$, we will compare it with its Kelvin 
transform $v_{\bp, \la}$ defined by (\ref{1.7}) for $\bp\in {\mathbb S}^n$ 
and $0<\la<\pi$. The first result indicates that the comparison 
is always possible if $\la>0$ is small and $v$ is regular at $\bp$.

\begin{Lemma}\label{L3.1} 
Assume that $n\ge 3$ and that $\Gamma$ is a closed subset of 
${\mathbb S}^n$. If $v\in C^1({\mathbb S}^n\backslash \Ga)$  and $v\ge 
c_0$ on ${\mathbb S}^n\backslash \Ga$ for some constant $c_0>0$, then 
for each $\bp\in {\mathbb S}^n\setminus\Gamma$ there exists 
$0<\la_\bp<\pi/2$ such that 
$$ 
v_{\bp, \la}\le v \quad\mbox{on } \Sigma_{\bp, 
\la}\setminus\Gamma \mbox{ for each } 0<\la<\la_\bp. 
$$ 
\end{Lemma} 
 
\begin{proof} 
Since $\bp\not\in \Gamma$ and $\Gamma$ is closed in ${\mathbb S}^n$, 
there exists $0<\lambda_0<\frac{\pi}{2}$ such that $\Gamma\subset 
\Sigma_{\bp, \lambda_0}$. Then for $0<r<\lambda_0$ we have 
\begin{align*} 
\frac{\p}{\p r}&\left\{(1+\cos^2\lambda-2\cos\lambda \cos 
r)^{\frac{n-2}{4}}v(r, \omega)\right\}\\ 
=&(1+\cos^2\lambda-2\cos\lambda\cos r)^{\frac{n-2}{4}}\left[ 
v_r+\frac{n-2}{2}\frac{\cos\lambda\sin 
r}{1+\cos^2\lambda-2\cos\lambda\cos r} v\right] 
\end{align*} 
Noting that $\sup_{B_{\la_0}(\bp)} |\nabla_{{\mathbb S}^n} v|$ is 
finite and $v\ge c_0$ on ${\mathbb S}^n\setminus \Gamma$,  there 
exists $0<\la_1<\lambda_0$ such that for $0<\lambda<r<\la_1$ 
$$ 
\frac{\p}{\p r}\left\{(1+\cos^2\lambda-2\cos\lambda \cos 
r)^{\frac{n-2}{4}}v(r, \omega)\right\}>0. 
$$ 
This implies that for $0<\lambda<r<\la_1$ 
$$ 
(1+\cos^2\lambda-2\cos\lambda\cos h_\lambda(r))^{\frac{n-2}{4}} 
v(h_\lambda(r), \omega)\le (1+\cos^2\lambda-2\cos\lambda\cos 
r)^{\frac{n-2}{4}}v(r, \omega). 
$$ 
Since 
$$ 
1+\cos^2\lambda-2\cos\lambda\cos 
h_\lambda(r)=\frac{\sin^4\lambda}{1+\cos^2\lambda-2\cos\lambda\cos 
r} 
$$ 
We therefore conclude that 
\begin{equation}\label{3.1} 
v_{\bp, \la}(r, \omega)\le v(r, \omega) \quad \mbox{if } 
0<\la<r<\la_1 \mbox{ and } \omega\in {\mathbb S}^{n-1}. 
\end{equation}

Next we can find a constant $C_1$ such that for $(r, \omega)\in 
\Sigma_{\bp, \la_1}\setminus \Ga$ 
\begin{align*} 
\frac{v_{\bp, \la}(r, \omega)}{v(r, \omega)}=& 
\left(\frac{\sin^2\la}{1+\cos^2\la-2\cos\la\cos 
r}\right)^{\frac{n-2}{2}} \frac{v(h_\la(r), \omega)}{v(r, \omega)}\\ 
\le &C_1 
\left(\frac{\sin^2\la}{1+\cos^2\la-2\cos\la\cos r} 
\right)^{\frac{n-2}{2}}, 
\end{align*} 
where we used the facts that $v\ge c_0>0$ and $v\circ \varphi_{\bp, 
\lambda}$ is bounded on $\Sigma_{\bp, \la_1}$. So there exists 
$0<\la_2<\la_1$ such that 
$$ 
v_{\bp, \la}\le v \quad \mbox{on  } \Sigma_{\bp, 
\la_1}\setminus \Gamma \mbox{ for each  } 0<\lambda<\lambda_2. 
$$ 
Combining this with (\ref{3.1}) we thus complete the proof. 
\end{proof}

\begin{Lemma}\label{L3.2} 
Let ${\mathcal O}$ be a domain in ${\mathbb S}^n$ and $\bq\in 
\overline{\mathcal O}$. If $v\in C^2({\mathcal 
O}\setminus\{\bq\})\cap C^0(\overline{\mathcal O}\backslash\{\bq\})$ 
is a non-negative function such that $-{\mathcal L}_{{\mathbb 
S}^n}v+Cv\ge 0$ in ${\mathcal O}\backslash \{\bq\}$ for some non-negative constant $C$ and $v\ge c_0$ on $\p 
{\mathcal O}\backslash\{\bq\}$ for some constant $c_0>0$, then 
$v>c_1$ on $\overline{\mathcal O}\backslash\{\bq\}$ for some 
constant $c_1>0$. 
\end{Lemma} 
 
\begin{proof} 
Using the stereographic projection with respect to $\bq$, the conclusion is a consequence of Lemma \ref{L4.1}. 
\end{proof} 
 
Fix a point ${\mathbf p}\in {\mathbb S}^n$ and let $\Ga\subset {\mathbb S}^n\setminus 
B_{\pi/2}({\mathbf p})$ be a set consisting of discrete points. 
Let $g:({\mathbb S}^n\setminus\Ga)\times (0, \infty)\to {\mathbb R}$ be a continuous function. We consider the equation 
\begin{equation}\label{3.2} 
-{\mathcal L}_{{\mathbb S}^n} v=g(\theta, v) \quad \mbox{and} \quad v>0 \quad \mbox{on } 
{\mathbb S}^n\setminus \Gamma. 
\end{equation} 
If $v\in C^2({\mathbb S}^n\setminus \Gamma)$ is a solution of (\ref{3.2}), 
we define $0<\bar{\lambda}_{\mathbf p}\le \pi$ by 
$$ 
\bar{\la}_{\mathbf p}:=\sup\{\la\in (0, \pi]: v_{{\mathbf p}, \mu}\le v \mbox{ 
in } \Sigma_{{\mathbf p}, \mu}\setminus \Gamma \mbox{ for each } 
0<\mu<\la\}. 
$$ 
Since $\Ga$ is discrete,  if we assume $\inf_{{\mathbb S}^n\setminus\Gamma}v>0$ then, by using Lemma \ref{L3.1}, 
$\bar{\lambda}_{\bp}$ is well-defined. The next result shows that $\la_{\mathbf 
p}\ge \pi/2$ if $g$ satisfies (g3) and the following conditions: 
 
\begin{enumerate} 
 
\item[(g7)$_{\bp}$] 
For each $s>0$, $0<\la<\pi/2$ and $\theta\in 
\Sigma_{{\mathbf p}, \la}$, $g(\theta, s)\ge g(\varphi_{{\mathbf p}, \la}(\theta), s)$. 
 
\item[(g8)$_\bp$] Either $g(\theta, s)>g(\varphi_{{\mathbf p}, \la}(\theta), s)$ for any $s>0$, $0<\la<\pi/2$ and $\theta\in \Sigma_{{\mathbf p}, \la}\setminus \Ga$,  or for $0<\la<\pi/2$ and $\theta\in\Sigma_{{\mathbf p},\la}$, the function $s\to s^{-(n+2)/(n-2)} g(\theta, s)$ is strictly deceasing. 
\end{enumerate} 
 
\begin{Lemma}\label{L3.3} 
For $\bp\in{\mathbb S}^n$, let $\Gamma\subset{\mathbb S}^n\setminus 
B_{\pi/2}({\mathbf p})$ be a discrete set. Assume $g$ is continuous on 
$\left({\mathbb S}^n\setminus\Gamma\right)\times(0,\infty)$ and $g(\cdot,s)$ 
is bounded below on ${\mathbb S}^n\setminus\Ga$ for each $s\in (0,\infty)$. 
Assume also that $g$ satisfy (g3), $(g7)_\bp$ and $(g8)_\bp$. 
If $v\in C^2({\mathbb S}^n\setminus\Gamma)$ is a solution of (\ref{3.2}) with 
$\inf_{{\mathbb S}^n\setminus \Gamma}v>0$, then 
$\bar{\lambda}_{\mathbf p}\ge \pi/2$. 
\end{Lemma} 
 
\begin{proof} 
In the following we will use the abbreviations 
$$ 
\varphi_\la:=\varphi_{{\mathbf p},\la}, \,\, v_\la:=v_{{\mathbf p}, \la}, \,\, 
B_{\la}:=B_{\la}({\mathbf p}), \,\, 
\Sigma_\la:=\Sigma_{{\mathbf p},\la}, \,\, \bar{\la}:=\bar{\la}_{\mathbf p}. 
$$ 
By Lemma \ref{L3.1}, $\bar\la$ is well defined. We argue by contradiction and assume $\bar{\la}<\pi/2$. 
>From the definition of $\bar{\la}$ we have 
\begin{equation}\label{3.3} 
v_\la\le v  \quad \mbox{ in } \Sigma_\la\setminus\Gamma 
\mbox{ for each } 0<\la\le \bar{\la}. 
\end{equation} 
Moreover, from (\ref{3.2}), (\ref{1.7}) and (\ref{1.8}) it 
follows that on $\Sigma_{\bar{\la}}\setminus\Gamma$ there hold 
\begin{equation}\label{3.4} 
-v^{-\frac{n+2}{n-2}}{\mathcal L}_{{\mathbb 
S}^n}v=v^{-\frac{n+2}{n-2}}g(\theta, v) 
\end{equation} 
and 
\begin{equation}\label{3.5} 
-v_{\bar{\la}}^{-\frac{n+2}{n-2}} {\mathcal L}_{{\mathbb S}^n} 
v_{\bar{\la}}=(v\circ \varphi_{\bar{\la}})^{-\frac{n+2}{n-2}} 
g(\varphi_{\bar{\la}}(\theta), v\circ \varphi_{\bar{\la}}). 
\end{equation} 
 
We define 
$$ 
{\mathcal O}:=\left\{\theta\in\Sigma_{\bar{\la}}\setminus\Gamma: 
v(\theta)<v\circ \varphi_{\bar{\la}}(\theta)\right\}. 
$$ 
Let $v_s=sv+(1-s)v_{\bar\la}$, for $0\le s\le 1$. Use the 
technique developed in the proof of \cite[Lemma 2.2]{LZ}, it follows from (g7)$_\bp$, (g3), 
(\ref{3.4}) and (\ref{3.5}) that 
\begin{align}\label{3.6} 
0&\ge -(v\circ \varphi_{\bar{\la}})^{-\frac{n+2}{n-2}}g(\theta, (v\circ \varphi_{\bar{\la}}))+(v\circ \varphi_{\bar{\la}})^{-\frac{n+2}{n-2}}g(\varphi_{\bar{\la}}(\theta), (v\circ \varphi_{\bar{\la}}))\nonumber\\ 
&\ge -v^{-\frac{n+2}{n-2}}g(\theta, v)+(v\circ \varphi_{\bar{\la}})^{-\frac{n+2}{n-2}} 
g(\varphi_{\bar{\la}}(\theta), v\circ \varphi_{\bar{\la}}) \nonumber\\ 
&=\int_0^1\frac{d}{ds}\left(v_s^{-\frac{n+2}{n-2}}{\mathcal L}_{{\mathbb S}^n}v_s\right)ds\nonumber\\ 
&=\left(\int_0^1 v_s^{-\frac{n+2}{n-2}}ds\right){\mathcal 
L}_{{\mathbb S}^n}(v-v_{\bar\la})-\frac{n+2}{n-2}\left(\int_0^1 
v_s^{-\frac{2n}{n-2}}{\mathcal L}_{{\mathbb 
S}^n}v_sds\right)(v-v_{\bar\la}) 
\end{align} 
in ${\mathcal O}$. 
 
We first claim that 
\begin{equation}\label{3.7} 
v-v_{\bar \la}>0 \quad \mbox{in } {\mathcal O}. 
\end{equation} 
In fact, By (\ref{3.3}), (\ref{3.6}) and the strong maximum 
principle, either (\ref{3.7}) holds or $v-v_{\bar\la}\equiv 0$ in $\mathcal O$. If ${\mathcal O}\ne \Sigma_{\bar\la}\setminus 
\Ga$, (\ref{3.7}) is true since $\p {\mathcal O}\cap 
(\Sigma_{\bar\la}\setminus \Ga)\ne \emptyset$ and 
$v-v_{\bar\la}>0$ on $\p{\mathcal O}\cap 
(\Sigma_{\bar\la}\setminus \Ga)$. So we may assume ${\mathcal 
O}=\Sigma_{\bar\la}\setminus\Ga$. If $v=v_{\bar\la}$ in 
$\Sigma_{\bar\la}\setminus \Ga$, then by (\ref{3.6}) we must have 
$$ 
v^{-(n+2)/(n-2)} g(\theta, v)\equiv(v\circ \varphi_{\bar \la})^{-(n+2)/(n-2)} g(\varphi_{\bar\la}(\theta), v\circ \varphi_{\bar\la}) \quad \mbox{in } \Sigma_{\bar\la}\setminus\Ga. 
$$ 
But (g8)$_\bp$ implies that this can not happen. We therefore obtain (\ref{3.7}).

We next claim that 
\begin{equation}\label{3.8} 
\left.\frac{\p}{\p r}(v-v_{\bar{\la}})\right|_{\p 
B_{\bar{\la}}}>0. 
\end{equation} 
>From (\ref{3.3}) we know the left hand side of (\ref{3.8}) is 
always non-negative. Suppose (\ref{3.8}) is not true, then there 
exists $\bp_0\in \p B_{\bar{\la}}$ such that $\frac{\p}{\p 
r}(v-v_{\bar{\la}})(\bp_0)=0$. However, direct calculation 
shows 
$$ 
\frac{\p}{\p r} (v\circ \varphi_{\bar{\la}}- v_{\bar{\la}}) 
=(n-2)(\cot \bar{\la}) v(\bp_0). 
$$ 
Since $0<\bar{\la}<\pi/2$, we then have 
$$ 
\frac{\p}{\p r}(v-v\circ \varphi_{\bar{\la}})(\bp_0)<0. 
$$ 
Since $v=v\circ \varphi_{{\bar \la}}$ on $\p B_{\bar{\la}}$, we 
conclude that ${\mathcal U}_{\bp_0}\cap 
\Sigma_{\bar{\la}}\subset{\mathcal O}$ for some neighborhood 
${\mathcal U}_{\bp_0}$ of $\bp_0$ in ${\mathbb S}^n$. Thus from 
(\ref{3.6}) and the Hopf Lemma it follows that $\frac{\p}{\p 
r}(v-v_{\bar{\la}})(\bp_0)>0$ which is a contradiction. We 
therefore obtain (\ref{3.8}). 
 
Noting that $v$ and $v_{\bar{\la}}$ are $C^2$ near $\p 
\Sigma_{\bar{\la}}$. By using (\ref{3.8}) it is easy to 
find $\bar{\la}<\la_0<\pi/2$ such that 
\begin{equation}\label{3.9} 
v_\la< v \quad\mbox{on }  B_{\la_0}\setminus 
\overline{B}_{\la}  \mbox{ for } \bar{\la}\le 
\la<\la_0. 
\end{equation} 
 
We still need to consider the points in $\Sigma_{\la_0} 
\setminus \Ga$. By using (\ref{1.7}) and the 
definition of ${\mathcal O}$, it is easy to see that there is a 
positive constant $\alpha_0$ such that 
\begin{equation}\label{3.10} 
v-v_{\bar{\la}}\ge \alpha_0 \quad \mbox{on } \Sigma_{\la_0} 
\setminus ({\mathcal O}\cup \Ga). 
\end{equation} 
This together with (\ref{3.7}) implies that 
\begin{equation}\label{3.11} 
v-v_{\bar{\la}}\ge 
c_0>0 \quad \mbox{on } \p ({\mathcal O}\cap \Sigma_{\la_0})\setminus \Ga 
\end{equation} 
for some constant $c_0>0$.

We observe that since $v$ is $C^2$ on the compact set ${\mathbb S}^n\setminus\Sigma_{\la_0}$, 
$$ 
v_{\bar\la}\le v\le v\circ \varphi_{\bar\la}\le C\quad\mbox{in } {\mathcal O}\cap \Sigma_{\la_0}. 
$$ 
Since $g(\cdot,C)$ is bounded below in ${\mathbb S}^n\setminus\Ga$, it follows from (g3) and (\ref{3.2}) that 
$$ 
-{\mathcal L}_{{\mathbb S}^n}v=g(\theta,v)\ge C^{-\frac{n+2}{n-2}}g(\theta,C)v^{\frac{n+2}{n-2}}\ge -Cv 
\quad\mbox{in } {\mathcal O}\cap\Sigma_{\la_0}. 
$$ 
It is easy to see  from (\ref{3.5}) that 
$$ 
-{\mathcal L}_{{\mathbb S}^n}v_{\bar\la}\ge -Cv_{\bar\la} \quad\mbox{in } 
{\mathcal O}\cap\Sigma_{\la_0}. 
$$ 
Note that $v$ and $v_{\bar\la}$ have positive lower and upper bounds in 
${\mathcal O}\cap\Sigma_{\la_0}$, we can use the condition 
$\inf_{{\mathbb S}^n\setminus\Ga}v>0$ and (\ref{3.6}) to obtain, 
for some positive constant $C$, 
$$ 
-{\mathcal L}_{{\mathbb S}^n}(v-v_{\bar\la})+C(v-v_{\bar\la})\ge 0,\quad\mbox{in } {\mathcal O}\cap\Sigma_{\la_0}. 
$$ 
Since $\Ga$ is discrete, by using (\ref{3.3}),  (\ref{3.11}) and  Lemma \ref{L3.2} we have 
$$ 
\inf_{{\mathcal O}\cap\Sigma_{\la_0}}\left(v-v_{\bar{\la}}\right)>0. 
$$ 
This together with (\ref{3.10}) implies, Thus for some constant $c_1>0$, 
$$ 
v-v_{\bar{\la}}\ge c_1 \quad \mbox{on } \Sigma_{\la_0} 
\setminus \Ga 
$$ 
Using (\ref{1.7}) then we can find 
$0<\varepsilon<\la_0-\bar{\la}$ such that 
\begin{equation}\label{3.12} 
v_{\la}<v \quad \mbox{on } \Sigma_{\la_0}\setminus \Ga 
\mbox{ for each } \bar{\la}\le 
\la<{\bar \la}+\epsilon. 
\end{equation} 
 
Combining (\ref{3.3}), (\ref{3.9}) and (\ref{3.12}) we have 
$$ 
v_{\la}\le v \mbox{ in } \Sigma_{\la}\setminus\Ga \mbox{ for each } 0<\la<\bar{\la}+\varepsilon. 
$$ 
This contradicts the definition of $\bar \la$. Hence 
$\bar{\la}\ge \pi/2$. 
\end{proof}

\subsection{Proof of Theorem \ref{T1.4}}

Now we are ready to give the proof of Theorem \ref{T1.4}. 
Under conditions (g1) and (g2), one can see that (g7)$_\bp$ and (g8)$_\bp$ are satisfied with $\bp={\mathbf s}$. 
So we may apply Lemma \ref{L3.3} to conclude that 
$$ 
v\ge v_{\mathbf s, \pi/2} \quad \mbox{on } \Sigma_{\mathbf s, \pi/2}\setminus\{{\mathbf n}\}. 
$$ 
Note that $\varphi_{\mathbf s, \pi/2}$ is a mirror reflection and 
$|J_{\varphi_{\mathbf s, \pi/2}}|=1$, we have from (\ref{3.2}), (g1), (g2) and (g4) that 
\begin{equation}\label{3.13} 
-{\mathcal L}_{{\mathbb S}^n} (v-v_{\mathbf s, \pi/2}) =g(\theta,v)-g(\varphi_{\mathbf s, \pi/2}(\theta), v_{\mathbf s, \pi/2})>0\quad  \mbox{in }  \Sigma_{\mathbf s, \pi/2}\setminus \{{\mathbf n}\}. 
\end{equation} 
Therefore by the strong maximum principle and the Hopf lemma we have 
\begin{equation}\label{3.14} 
v>v_{\mathbf s ,\pi/2} \quad \mbox{in } \Sigma_{\mathbf s,\pi/2}\setminus\{{\mathbf n}\} 
\end{equation} 
and 
\begin{equation}\label{3.15} 
\left.\frac{\p}{\p r}(v-v_{\mathbf s,\pi/2}) \right|_{\p B_{\pi/2}(\mathbf s)}>0. 
\end{equation} 
By Lemma \ref{L3.2}, 
\begin{equation}\label{3.16} 
\inf_{\theta\to\bn}(v(\theta)-v_{\mathbf s,\pi/2}(\theta))>c, 
\end{equation}for some constant $c>0$. 
 
For any $\bq\in\p B_{\pi/2}(\mathbf s)$, $\bq=(\pi/2,\omega_0)$ in the geodesic polar coordinate with respect to $\mathbf s$. 
Let $\bp_t=(t\pi/2,\omega_0)$ for $0\le t\le 2$. 
 
{\it Claim 1.} There exists $\varepsilon>0$, such that for any $0\le t<\varepsilon$, there holds 
\begin{equation}\label{3.17} 
v\ge v_{\mathbf p_t ,\pi/2} \quad \mbox{in } \Sigma_{\mathbf p_t ,\pi/2}\setminus\{{\mathbf n}\}. 
\end{equation} 
The claim follows easily from (\ref{3.14}), (\ref{3.15}) and 
(\ref{3.16}). For readers' convenience, we include a proof by 
contradiction argument. If it is not true, then there exists a 
sequence $0<t_i\to 0$ and a sequence  $\{\theta_i\}$ with 
$\theta_i \in\Sigma_{{\mathbf p}_{t_i},\pi/2}\setminus\{{\mathbf 
n}\}$ such that 
\begin{equation}\label{3.18} 
v(\theta_i)<v_{{\mathbf p}_{t_i},\pi/2}(\theta_i). 
\end{equation} 
By taking a subsequence if necessary, we may assume $\{\theta_i\}$ 
converges to some point $\theta_0\in\overline{\Sigma_{\mathbf 
s,\pi/2}}$. If $\theta_0\in\Sigma_{\mathbf s, 
\pi/2}\setminus\{{\mathbf n}\}$, then $v(\theta_0)\le v_{\mathbf 
s,\pi/2}(\theta_0)$ which violates (\ref{3.14}). If 
$\theta_0={\mathbf n}$, (\ref{3.16}) implies that $(v-v_{{\mathbf 
p}_i, \pi/2})(\theta_i)>c/2$ for large $i$,  a contradiction to 
(\ref{3.18}). Finally we assume $\theta_0\in \p B_{\pi/2}(\mathbf 
s)$. Let $\bar{\theta}_i$ be the closest point on 
$\p\Sigma_{{\mathbf p}_{t_i}, \pi/2}$ to $\theta_i$. Note that $v(\bar 
\theta_i) =v_{{\mathbf p}_{t_i}, \pi/2}(\bar \theta_i)$, by using 
(\ref{3.18}) we have 
$$ 
\frac{\p}{\p r_i} (v-v_{{\mathbf p}_{t_i}, \pi/2})(\tilde{\theta}_i)<0 
$$ 
for some $\tilde{\theta}_i$ between $\theta_i$ and $\bar \theta_i$ 
on the geodesic line connecting $\theta_i$ and $\bar \theta_i$, 
where $r_i$ denotes the geodesic distance from ${\mathbf p}_{t_i}$. Note 
that $\tilde{\theta}_i \rightarrow \theta_0$. We obtain 
$$ 
\frac{\p}{\p r}(v-v_{\mathbf s, \pi/2})(\theta_0)\le 0. 
$$ 
which violates (\ref{3.15}). We thus prove the claim. 
 
We define 
$$ 
\bar t: =\sup\left\{t\in (0, 1): v\ge v_{{\mathbf p}_\tau, \pi/2} 
\mbox{ in } \Sigma_{{\mathbf p}_\tau, \pi/2}\setminus\{{\mathbf n}\} 
\mbox{ for all } 0\le\tau<t\right\}. 
$$ 
>From Claim 1, $\bar t$ is well defined and $\bar t>0$. 
 
{\it Claim 2}. $\bar t= 1$. 
 
Suppose $\bar t<1$, by continuity of $v$, we have 
$$ 
v\ge v_{{\mathbf p}_{\bar t},  \pi/2} \quad \mbox{in } \Sigma_{{\mathbf p}_{\bar t},\pi/2}\setminus\{{\mathbf n\}}. 
$$ 
the conditions on $g$ imply 
\begin{equation}\label{3.19} 
-{\mathcal L} _{{\mathbb S}^n} (v-v_{{\mathbf p}_{\bar t},\pi/2})> 0 \quad \mbox{in } \Sigma_{{\mathbf p}_{\bar t}, \pi/2}\setminus\{{\mathbf n}\} 
\end{equation} 
and 
$$ 
v\neq v_{{\mathbf p}_{\bar t}, \pi/2}\quad \mbox{in } \Sigma_{{\mathbf p}_{\bar t}, \pi/2}\setminus\{{\mathbf n}\}. 
$$ 
Similar to the proof of Claim 1, there exists $\varepsilon_{\bar t}>0$ such that 
for all $\bar t\le \mu<\bar t+\varepsilon_{\bar t}<1$, 
$$ 
v\ge v_{{\mathbf p}_\mu, \pi/2} \quad \mbox{in } \Sigma_{{\mathbf p}_\mu,\pi/2}\setminus\{{\mathbf n\}}, 
$$ 
which contradicts the definition of $\bar t$. Claim 2 thus follows. 
 
By continuity of $v$, we finally obtain $v\ge v_{{\mathbf q}, \pi/2}$ in 
$\Sigma_{{\mathbf q}, \pi/2}$ for ${\mathbf q}\in \p 
B_{\pi/2}({\mathbf s})$. Since $\bq$ is arbitrarily chosen on $\p B_{\pi/2}(\mathbf s)$, the proof is complete.

\subsection{Proof of Theorem \ref{T1.5} and its corollaries} 
 
We first use Lemma \ref{L3.3} to give the proof of Theorem 
\ref{T1.5}. 
 
\begin{proof}[Proof of Theorem \ref{T1.5}] 
Using the conditions (g1), (g5) and (g6), it is easy to see 
(g7)$_\bp$ and (g8)$_\bp$ are satisfied for every $\bp \in \p 
B_{\pi/2}({\mathbf s})$. Therefore Lemma \ref{L3.3} with 
$\Ga=\{\bn, {\mathbf s}\}$ implies that ${\bar \lambda}_\bp\ge 
\pi/2$ for each $\bp \in \p B_{\pi/2}({\mathbf s})$. This in 
particular implies $v_{\bp, \pi/2}\le v$ on $\Sigma_{\bp, 
\pi/2}\setminus\{\bn, {\mathbf s}\}$ for each $\bp\in \p 
B_{\pi/2}({\mathbf s})$. Consequently $v$ is rotationally symmetric about 
the line through $\bn$ and ${\mathbf s}$. 
\end{proof}

In order to use Theorem \ref{T1.5} to prove Corollary \ref{C1.3}, 
we will show $f(v)>0$ on ${\mathbb S}^n$ for any solution $v\in 
C^2({\mathbb S}^n)$ of (\ref{1.13}), thus condition (\ref{1.10}) is 
satisfied due to Lemma \ref{L3.2}. This is given by the following 
simple observation. 
 
\begin{Lemma}\label{L3.4} 
Let $f$ satisfy the conditions given in Corollary \ref{C1.3}. 
If $v\in C^2({\mathbb S}^n)$ is a solution of 
(\ref{1.13}), then $f(v)>0$ on ${\mathbb S}^n$. 
\end{Lemma} 
 
\begin{proof} 
Suppose it is not true, then there is $\bar{\theta}\in {\mathbb S}^n$ 
such that $f(v(\bar{\theta}))\le 0$. Let $\hat{\theta}\in {\mathbb 
S}^n$ be a point such that $v(\hat{\theta})=\max_{{\mathbb S}^n} v$. 
Then it follows from the condition on $f$ that 
$$ 
v(\hat{\theta})^{-(n+2)/(n-2)}f(v(\hat{\theta}))\le v(\bar{\theta})^{-(n+2)/(n-2)} 
f(v(\bar{\theta}))\le 0. 
$$ 
On the other hand, by the 
maximality of $v$ at $\hat{\theta}$ we have $\Delta 
v(\hat{\theta})\le 0$, it follows from (\ref{1.13}) that 
$$ 
f(v(\hat{\theta}))\ge \frac{n(n-2)}{4} v(\hat{\theta})>0. 
$$ 
We thus derive a contradiction. 
\end{proof}

\begin{proof}[Proof of Corollary \ref{C1.3}] 
It is clear that $g(\theta, s):=f(s)$ satisfies (g1), (g3), (g5) and (g6). 
Since $v\in C^2({\mathbb S}^n)$ and Lemma \ref{L3.4} implies $f(v)>0$ 
on ${\mathbb S}^n$, we can use Theorem \ref{T1.5} to conclude $v$ is 
rotationally symmetric about the line through $\bp$ and $-\bp$ for any $\bp\in 
{\mathbb S}^n$. Therefore $v$ must be constant. 
\end{proof}

\begin{proof}[Proof of Corollary \ref{C1.4}] 
(i) When $\beta<n(n-2)/4$, the function $g(\theta,s)$ defined by 
(\ref{1.14}) satisfies (g1), (g3), (g5) and (g6) 
and is positive for $s\in (0, \infty)$. 
Thus we can apply Theorem \ref{T1.5} to conclude that any 
smooth solution $v\in C^2({\mathbb S}^n\backslash\{\bn, {\mathbf s}\})$ of 
(\ref{1.15}) is rotationally symmetric about the line through $\bn$ and 
${\mathbf s}$. 
 
Given a solution $v\in C^2({\mathbb S}^n\backslash \{\bn, {\mathbf s}\})$ of 
(\ref{1.13}), the function 
$$ 
u(x):=\xi(x) v(\pi_{\bn}^{-1}(x)) 
$$ 
is then radially symmetric in ${\mathbb R}^n\backslash\{0\}$. Write $u(x)=u(r)$ with 
$r=|x|$, then $u$ satisfies the ordinary differential equation 
\begin{equation}\label{3.20} 
u''+\frac{n-1}{r}u'+\frac{n(n-2)-4\beta}{(1+r^2)^2}u+ 
u^{\frac{n+2}{n-2}}=0,\,\, u(r)>0 \mbox{ for } 0<r<\infty. 
\end{equation}

(ii) Now we assume $\beta\le 0$ and  (\ref{1.13}) has a solution 
$v\in C^2({\mathbb S}^n\backslash\{\bn, {\mathbf s}\})$. This implies that 
(\ref{3.20}) has a solution $u$. Let 
$\varphi(r):=(1+r^2)^{\frac{n-2}{2}}u(r)$, then 
\begin{equation}\label{3.21} 
\left(\frac{r^{n-1}}{(1+r^2)^{n-2}}\varphi'\right)'= 
-\frac{r^{n-1}}{(1+r^2)^n}\left(\varphi^{\frac{n+2}{n-2}}-4\beta 
\varphi\right) \quad \mbox{ for } 0<r<\infty. 
\end{equation} 
Using $\beta\le 0$ we have for any $r>\varepsilon>0$ that 
\begin{equation}\label{3.22} 
\frac{r^{n-1}}{(1+r^2)^{n-2}}\varphi'(r)<\frac{\varepsilon^{n-1}}{(1+ 
\varepsilon^2)^{n-2}}\varphi'(\varepsilon). 
\end{equation} 
If 
$$ 
\lim_{r\rightarrow 0} \frac{r^{n-1}}{(1+r^2)^{n-2}}\varphi'(r)>0, 
$$ 
then there exists a number $\alpha_0>0$ such that $\varphi'(r)\ge 
\alpha_0 r^{1-n}$ for small $r>0$. Therefore for any small 
$r>\varepsilon>0$ there holds 
$$ 
\varphi(r)\ge \varphi(r)-\varphi(\varepsilon)\ge 
\frac{\alpha_0}{n-2}\left(\varepsilon^{2-n}-r^{2-n}\right). 
$$ 
Fixing $r$ and letting $\varepsilon\rightarrow 0$ we then derive a 
contradiction. Therefore 
$$ 
\lim_{r\rightarrow 0} \frac{r^{n-1}}{(1+r^2)^{n-2}}\varphi'(r)\le 
0. 
$$ 
It then follows from (\ref{3.22}) that $\varphi'(r)<0$ for all 
$r>0$. Let $\tilde{\varphi}(r):=\varphi(\frac{1}{r})$, then 
$\tilde{\varphi}'(r)>0$ for all $r>0$. However, by direct 
calculation one can see that $\tilde{\varphi}$ is also a solution 
of (\ref{3.21}). Therefore the above argument applies to 
$\tilde{\varphi}$ and shows that $\tilde{\varphi}'(r)<0$ for all 
$r>0$. We thus derive a contradiction.

(iii) In order to show that (\ref{1.15}) has infinitely many 
solutions, it is equivalent to showing that (\ref{3.20}) has 
infinitely many solutions defined on $(0, \infty)$. To this end, 
for any function $u$ defined on $(0, \infty)$ we define 
$$ 
w(t)=e^{-\frac{n-2}{2}t} u(e^{-t}), \quad t\in (-\infty, \infty). 
$$ 
By an easy calculation one can see that $u$ satisfies (\ref{3.20}) 
if and only if $w$ satisfies 
\begin{equation}\label{3.23} 
w''-\left(\frac{n-2}{2}\right)^2 w+w^{\frac{n+2}{n-2}}+ 
\frac{c_\beta e^{-2t}}{\left(1+e^{-2t}\right)^2} w=0, \quad w>0 
\mbox{ on } (-\infty, \infty), 
\end{equation} 
where $c_\beta:=n(n-2)-4\beta$. Therefore, it suffices to show 
that (\ref{3.23}) has infinitely many positive solutions. 
 
Let us introduce a function $h(\cdot, \cdot): {\mathbb R}_+\times 
{\mathbb R}\to {\mathbb R}$ by 
$$ 
h(a, b):=b^2-\left(\frac{n-2}{2}\right)^2 a^2+\frac{n-2}{n} 
a^{\frac{2n}{n-2}}, \quad \forall (a, b)\in {\mathbb R}_+\times 
{\mathbb R}. 
$$ 
Since $\beta>\frac{n-2}{2}$, there are infinitely many $(a, b)\in 
{\mathbb R}^+\times {\mathbb R}$ such that 
\begin{equation}\label{3.24} 
h(a, b)+\frac{1}{4} c_\beta a^2<0. 
\end{equation} 
Let us fix one of them and consider the initial value problem 
$$ 
\left\{\begin{array}{lll} w''-\left(\frac{n-2}{2}\right)^2 
w+w^{\frac{n+2}{n-2}}+c_\beta 
\frac{e^{-2t}}{\left(1+e^{-2t}\right)^2} w=0\\ 
\\ 
w(0)=a,\quad w'(0)=b. 
\end{array}\right. 
$$ 
By the local existence theory for ordinary differential equations, 
it has a unique solution $w$ defined on some interval containing 
$t=0$. Let $(-B, A)$ be the largest interval on which $w$ exists 
and $w(t)>0$. Since $w(0)=a>0$, $A$ and $B$ must be positive. It 
remains only to show that $A=\infty$ and $B=\infty$. In the 
following we will only prove $A=\infty$, since $B=\infty$ can be 
proven in the same way. 
 
Suppose $A<\infty$ and consider the function 
$$ 
h(t):=h(w(t), w'(t)) \quad \mbox{on } [0, A). 
$$ 
>From the equation satisfied by $w$ it follows that 
$$ 
h'(t)=-c_\beta \frac{e^{-2t}}{\left(1+e^{-2t}\right)^2} 
\left[w(t)^2\right]'. 
$$ 
Therefore, for $t\in [0, A)$ one has, by integration by parts and 
noting that $\frac{e^{-2t}}{\left(1+e^{-2t}\right)^2}$ is 
non-increasing on $[0, \infty)$, 
\begin{align*} 
h(t)&-h(0)\\ 
=&-c_\beta\left\{\frac{e^{-2t}}{\left(1+e^{-2t}\right)^2} 
w(t)^2-\frac{1}{4} w(0)^2-\int_0^t w(s)^2 
\left[\frac{e^{-2s}}{\left(1+e^{-2s}\right)^2}\right]'ds 
\right\}\\ 
\le& \frac{1}{4} c_\beta w(0)^2. 
\end{align*} 
Consequently, by using (\ref{3.24}), 
\begin{equation}\label{3.25} 
h(t)\le h(a, b)+\frac{1}{4} c_\beta a^2<0 \quad \mbox{on } [0, A). 
\end{equation} 
This together with the equation of $w$ implies that 
$$ 
w+|w'|+|w''|\le C \quad \mbox{on } [0, A) 
$$ 
for some positive constants $C$. Therefore $w(A)$ and $w'(A)$ are 
well-defined and hence $w$ has definition on a larger interval 
$[0, A+\varepsilon)$ for some $\varepsilon>0$. But from 
(\ref{3.25}) one can see that $w(A)>0$. Thus, by continuity, $w>0$ 
on $[A, A+\varepsilon)$ for some smaller $\varepsilon>0$. This 
contradicts the definition of $A$. Therefore $A=\infty$. 
\end{proof}

\subsection{Proof of Theorem \ref{T1.8}, Theorem \ref{T1.9} and Theorem \ref{T1.10}} 
 
\begin{proof}[Proof of Theorem \ref{T1.8}] 
Let ${\mathbf s}$ be the south pole 
of ${\mathbb S}^n$. As before we define 
$$ 
\bar{\la}_{\mathbf s}:=\sup\left\{\la\in (0, \pi): v_{{\mathbf s}, \mu}\le v \mbox{ in } 
\Sigma_{{\mathbf s}, \mu}\setminus\{\bn\} \mbox{ for  each } 
0<\mu<\la\right\}. 
$$ 
Let $(r,\omega)$ be the geodesic polar coordinates on ${\mathbb 
S}^n$ with respect to $\bf s$. Then the conditions on $K$ are equivalent to saying that 
$K$ is non-constant on ${\mathbb S}^n\setminus\{\bn\}$ and 
for each fixed $\omega\in {\mathbb S}^{n-1}$ the function $r\to K(r,\omega)$ is 
non-decreasing for $0<r<\pi$. 
 
By condition (\ref{1.23}), it follows from Lemma \ref{L3.1} that 
$\bar{\la}:= \bar{\la}_{\mathbf s}>0$ is well defined. We claim 
that $\bar{\la}=\pi$. If $\bar{\la}<\pi$, then from 
(\ref{1.22}), (\ref{1.7}) and (\ref{1.8}) that 
$$ 
-v^{-\frac{n+2}{n-2}}{\mathcal L}_{{\mathbb 
S}^n}v(r,\omega)+v_{{\mathbf 
s},\bar{\la}}^{-\frac{n+2}{n-2}}{\mathcal L}_{{\mathbb 
S}^n}v_{{\mathbf s}, \bar{\la}}(r,\omega) 
=K(r,\omega)-K(h_{\bar{\la}}(r),\omega). 
$$ 
Following the proof of Lemma \ref{L3.3} and using the conditions on 
$K$ we have 
$$ 
0\ge \left(\int_0^1 v_s^{-\frac{n+2}{n-2}}ds\right){\mathcal 
L}_{{\mathbb S}^n}(v-v_{{\mathbf s}, \bar\la})-\frac{n+2}{n-2}\left(\int_0^1 
v_s^{-\frac{2n}{n-2}}{\mathcal L}_{{\mathbb 
S}^n}v_sds\right)(v-v_{{\mathbf s},\bar\la}) 
$$ 
in $\Sigma_{{\mathbf s}, \bar{\la}}\setminus\{\bn\}$. 
 
Note that $v-v_{{\mathbf s}, \bar{\la}}\ge  0$ on 
$\Sigma_{{\mathbf s}, \bar{\la}}\setminus\{\bn\}$. 
Since $K$ is non-constant on ${\mathbb S}^n\setminus\{\bn\}$, 
we must have $v-v_{{\mathbf s}, 
\bar{\la}}\neq 0$ on $\Sigma_{{\mathbf s}, 
\bar{\la}}\setminus\{\bn\}$. 
Thus it follows from the strong maximum principle and the Hopf lemma that 
$$ 
v-v_{{\mathbf s}, \bar{\la}}> 0 \quad \mbox{on }  \Sigma_{{\mathbf s}, \bar{\la}}\setminus\{\bn\} 
$$ 
and 
$$ 
\frac{\p}{\p r}(v-v_{{\mathbf s}, \bar{\la}})>0 \quad \mbox{on } \p \Sigma_{{\mathbf s}, \bar{\la}}. 
$$ 
 
Next we show that 
$$ 
\inf_{\theta\to {\mathbf n}}\left(v-v_{{\mathbf s},\bar{\la}}\right)>0. 
$$ 
Choose $0<r_0$ small such that $\bar{\la}+r_0<\pi$. Define 
$$ 
{\mathcal O}:=\left\{\theta\in B_{r_0}({\mathbf n})\setminus\{{\mathbf n}\}: 
v(\theta)<2v_{{\mathbf s},\bar{\la}}(\theta)\right\}. 
$$ 
By the fact that $v$ has positive lower bound, we can obtain $v-v_{{\mathbf s},\bar{\la}}\ge c_1$ on 
$\overline{B_{r_0}({\mathbf n})}\setminus\left({\mathcal O}\cup\{{\mathbf n}\}\right)$ for some positive constant $c_1$; 
moreover, $v$ and $v_{{\mathbf s},\bar{\la}}$ have lower and upper bounds in ${\mathcal O}$. 
Following the proof of Lemma \ref{L3.3}, we can obtain that for some large positive constant $C$ 
$$ 
-{\mathcal L}_{{\mathbb S}^n}(v-v_{{\mathbf s},\bar{\la}})+C(v-v_{{\mathbf s},\bar{\la}})\ge 0 
\quad\mbox{in}\quad {\mathcal O}. 
$$ 
Lemma \ref{L3.2} implies that $v-v_{{\mathbf s},\bar{\la}}>c_2$ in ${\mathcal O}$ 
for some positive constant $c_2$. We thus obtain 
$$ 
v-v_{{\mathbf s},\bar{\la}}\ge\min(c_1,c_2)>0 \quad\mbox{in  }  B_{r_0}({\mathbf n})\setminus\{{\mathbf n}\}. 
$$ 
 
Now we can argue as in the proof of Lemma \ref{L3.3} to show that there exists 
$\varepsilon>0$ such that $v\ge v_{{\mathbf s}, \la}$ on $\Sigma_{{\mathbf s}, \la}\setminus\{\bn\}$ for each $0<\la<\bar{\la}+\varepsilon$. This contradicts the definition of $\bar{\la}$. Therefore $\bar{\la}_{\mathbf s}=\pi$.

Next we are going to show (\ref{1.24}). In the above we have shown that 
\begin{equation}\label{3.26} 
v\ge v_{{\mathbf s}, \la} \quad \mbox{on } \Sigma_{{\mathbf s}, \la}\setminus\{\bn\} 
\quad \mbox{for each } 0<\la<\pi. 
\end{equation} 
For $\theta \in B_{\pi/2}(\bn)\setminus\{\bn\}$, let $(r, \omega)$ be its geodesic polar coordinate with respect to the south pole ${\mathbf s}$. Then we take $\pi/2<\la< \pi$ such that 
\begin{equation}\label{3.27} 
2\cos\la=(1+\cos^2\la)\cos r. 
\end{equation} 
This implies that $\varphi_{{\mathbf s},\la}(\theta)\in \p B_{\pi/2}(\bn)$. Moreover $\la\rightarrow \pi$ as $\theta\rightarrow \bn$. Let $C_1:=\min_{\p B_{\pi/2}(\bn)} v$, we have from (\ref{3.26}) and (\ref{3.27}) that 
$$ 
\liminf_{\theta\rightarrow \bn} d(\theta, \bn)^{(n-2)/2} v(\theta) 
\ge C_1 \liminf_{r\rightarrow \pi} \left(\frac{(\pi-r)\sin^2\la}{1+\cos^2\la-2\cos\la \cos r}\right)^{(n-2)/2}=C_1. 
$$ 
The proof is complete. 
\end{proof}

The proofs of Theorem \ref{T1.9} and Theorem \ref{T1.10} are based on the method used by  Li and Li in \cite{LL}. 
 
\begin{proof}[Proof of Theorem \ref{T1.9}] 
Since $0<v\in C^2({{\mathbb S}^n})$, by Lemma \ref{L3.1} the moving sphere 
procedure can start from south pole $\bf s$. So there exists 
$0<\lambda_1<\pi$, such that 
$$ 
v\ge v_{\bf s,\lambda}\quad\mbox{in } \Sigma_{\bf 
s,\lambda} \mbox{ for each } 0<\lambda<\lambda_1. 
$$ 
 
Let 
$$ 
\bar{\lambda}=\sup\{\lambda\in (0, \pi): v\ge v_{\bf 
s,\mu} \mbox{ in } \Sigma_{\bf s,\mu} \mbox{ for each } 0<\mu<\lambda\}. 
$$ 
We will show that $\bar\lambda=\pi$. If not, say $0<\bar\lambda<\pi$. 
Then by the conformal invariance and the condition on $K$ we have 
$$ 
f\left(\la(A_{v^{\frac{4}{n-2}}g_0})\right)-f\left(\la(A_{v_{\mathbf 
s,\bar\la}^{\frac{4}{n-2}}g_0})\right)=K(r,\omega)-K(h_{\bar{\la}}(r),\omega)\ge 
0 
$$ 
with strict inequality somewhere in $\Sigma_{\bf s,\bar\la}$. 
 
By the argument in \cite[Lemma 2.1]{LL}, there exists an elliptic 
operator $L$ such that 
$$ 
L(v-v_{\mathbf s,\bar\la})\ge 0 \quad\mbox{in}\quad \Sigma_{\bf 
s,\bar\la} 
$$ 
with strict inequality somewhere in $\Sigma_{\bf 
s,\bar\la}$. Noticing that $v\ge v_{\bf s,\bar\lambda}$ in 
$\Sigma_{\bf s,\bar\lambda}$, we can obtain that 
\begin{equation}\label{3.28} 
v-v_{\bf s,\bar\lambda}>0 \quad \mbox{in } \Sigma_{\bf s,\bar\lambda} 
\end{equation} 
and 
\begin{equation}\label{3.29} 
\frac{\partial (v-v_{\bf s,\bar\lambda})}{\partial r}>0 \quad \mbox{on } \p \Sigma_{\bf s, \bar \la}. 
\end{equation} 
 
Once (\ref{3.28}) and (\ref{3.29}) are established, we can argue as in 
the proof of Lemma \ref{L3.3} to show that we can move spheres 
beyond $\bar\lambda$, which contradict the definition of 
$\bar\lambda$. So $\bar\lambda=\pi$. Following the argument of Theorem \ref{T1.8}, 
we can see that $v$ blows up at ${\mathbf n}$. Therefore 
(\ref{1.31}) has no positive $C^2({\mathbb S}^n)$ solutions. 
\end{proof} 
 
\begin{proof}[Proof of Theorem \ref{T1.10}] 
Let $(r,\omega)$ be the geodesic polar coordinate with respect to 
$\bf s$. Let $\Sigma_{{\bf s}, \la}^+={\mathbb S}^n_+\cap \Sigma_{{\bf s}, \la}$. 
We define 
$$ 
\bar\lambda=\sup\{\lambda\in (0, \pi): v\ge v_{\bf s, \mu} \mbox{ in 
}\Sigma_{\bf s, \mu}^+  \mbox{ for all } 0<\mu<\lambda\}. 
$$ 
Since $v$ has positive lower bound on ${\mathbb 
S}^n_+$, it follows from Lemma \ref{L3.1} that $\bar{\lambda}$ is well-defined and $\bar{\lambda}>0$. 
We will show that $\bar\lambda=\pi$. If not, say $0<\bar\lambda<\pi$, then 
$$ 
v\ge v_{\bar\lambda} \quad\mbox{in}\quad \Sigma_{{\bf s}, \bar\lambda}^+. 
$$ 
Moreover, by the conformal invariance we have 
$$ 
\left\{ 
\begin{array}{l} 
f\left((A_{v_{{\bf s}, \bar\lambda}^{\frac{4}{n-2}}g_0})\right) 
=K(h_{\bar\la}(r),\omega)\quad \mbox{in } {\mathbb S}^n_+, \\ 
\\ 
\frac{\partial v_{{\bf s}, \bar\lambda}}{\partial\nu}=H(h_{\bar\la}(r),\omega) 
v_{{\bf s}, \bar\la}^{\frac{n}{n-2}}\quad \mbox{on } \partial {\mathbb S}^n_+. 
\end{array}\right. 
$$ 
By the argument in \cite[Lemma 2.1]{LL} and the condition $\nabla_{\frac{\p}{\p \theta_{n+1}}} K\ge 0$ 
on ${\mathbb S}^n_+$, there exists a linear elliptic operator $L$ such that 
$$ 
L\left(v-v_{{\bf s}, \bar{\lambda}}\right)\ge 0 \quad \mbox{in } \Sigma_{{\bf s}, \bar{\lambda}}^+. 
$$ 
Moreover 
$$ 
\frac{\p (v-v_{{\bf s}, \bar{\lambda}})}{\p \nu}=H(r, \omega) v^{\frac{n}{n-2}}-H(h_{\bar{\lambda}}(r), \omega) 
v_{{\bf s}, \bar{\lambda}}^{\frac{n}{n-2}} \quad \mbox{on } \p {\mathbb S}^n_+\cap \Sigma_{{\bf s}, \bar{\lambda}}^+. 
$$ 
Since $K$ or $H$ is non-constant and $\nabla_{\frac{\p}{\p \theta_{n+1}}} H\ge 0$ on $\p {\mathbb S}^n_+$, 
it follows from the strong maximum principle and the Hopf lemma that 
\begin{equation}\label{3.30} 
v>v_{{\bf s}, \bar\lambda}\quad \mbox{in } \Sigma_{{\bf s}, \bar\lambda}^+. 
\end{equation} 
Using the Hopf lemma again and \cite[Lemma 10.1]{LZ} we can also obtain 
$$ 
\frac{\partial(v-v_{{\bf s}, \bar\lambda})}{\partial 
r}|>0 \quad \mbox{on } \p \Sigma_{{\bf s}, \bar{\lambda}}\cap {\mathbb S}^n_+ 
$$ 
Thus there exists $\bar{\lambda}<\lambda_0<\pi$  such that 
$$ 
v\ge v_{{\bf s}, \lambda} \quad\mbox{on}\quad 
\Sigma_{{\bf s}, \lambda_0}^+\setminus \Sigma_{{\bf s}, \lambda}^+ \quad\mbox{for } 
\bar\lambda<\lambda<\lambda_0. 
$$ 
By using (\ref{3.30}) and the definition of $v_{{\bf s}, \lambda}$, we can find 
$0<\epsilon<\lambda_0-\bar\lambda$ such that 
$$ 
v\ge v_{{\bf s}, \lambda} \quad\mbox{on } \Sigma_{{\bf s}, \lambda_0}^+ \mbox{ for } 
\bar\lambda<\lambda<\bar\lambda+\epsilon. 
$$ 
Therefore 
$$ 
v\ge v_{{\bf s}, \lambda} \quad\mbox{on }\Sigma_{{\bf s}, \lambda}^+ \mbox{ for }\bar\lambda 
<\lambda<\bar\lambda+\epsilon 
$$ 
We thus derive a contradiction to the definition of $\bar\lambda$. Hence 
$\bar\lambda=\pi$. Consequently (\ref{1.32}) has no positive $C^2({\mathbb S}^n_+)$ 
solutions. 
\end{proof}

\subsection{\bf Proofs of Theorem \ref{T1.6} and Theorem \ref{T1.7}} 
 
Given a function $v$ on ${\mathbb S}^2$, we define, for each fixed 
$\bp\in {\mathbb S}^2$ and $0<\lambda<\pi$,  its Kelvin transform as 
$$ 
 v_{\bp, \lambda}=v\circ \varphi_{\bp, \lambda}+\frac{1}{2}\log 
|J_{\varphi_{\bp, \lambda}}|. 
$$ 
By conformal invariance it is known that 
\begin{equation}\label{3.31} 
-\Delta_{{\mathbb S}^2}v_{\bp, \lambda}+1=|J_{\varphi_{\bp, 
\lambda}}|\left((-\Delta_{{\mathbb S}^n}v)\circ \varphi_{\bp, 
\lambda}+1\right). 
\end{equation} 
 
The proof of Theorem \ref{T1.6} follows essentially the same idea 
as in the proof of Theorem \ref{T1.5}. We need to compare the 
functions $v$ and $v_{\bp, \lambda}$. Similar to Lemma \ref{L3.1} we 
have 
 
\begin{Lemma}\label{L3.5} 
Let $\Gamma$ be a closed subset of ${\mathbb S}^2$. If $v\in 
C^1({\mathbb S}^2\backslash \Gamma)$ and $v\ge -C_0$ on ${\mathbb 
S}^2\backslash \Gamma$ for some constant $C_0>0$, then for each 
$\bp\in {\mathbb S}^2\setminus \Gamma$ there exists 
$0<\lambda_{\bp}<\frac{\pi}{2}$ such that 
$$ 
v_{\bp, \lambda}\le v,\quad\mbox{on } \Sigma_{\bp, \lambda}\setminus 
\Gamma \mbox{ for each } 0<\lambda<\lambda_{\bp}. 
$$ 
\end{Lemma} 
 
The next two lemmas, similar to Lemma \ref{L3.2}, are used to deal with singularities. 
 
\begin{Lemma}\label{L3.6} 
Let $\mathcal O$ be an open set in ${\mathbb S}^2$ and $\bq\in\mathcal 
O$. If $v\in C^2({\mathcal O}\setminus\{\bq\})$ satisfies 
$-\Delta_{{\mathbb S}^2}v +1 >0$ in ${\mathcal O}\setminus\{\bq\}$ and 
\begin{equation}\label{3.32} 
\limsup_{\theta\in {\mathcal O}, \theta\rightarrow 
\bq}\frac{v(\theta)}{\log d(\theta, \bq)}\le 0, 
\end{equation} 
then $v>-C_0$ in ${\mathcal O}\setminus\{\bq\}$ for some constant 
$C_0>0$. 
\end{Lemma} 
\begin{proof} 
Fix an open set ${\mathcal O}'$ such that $\bq\in {\mathcal 
O}'\subset \overline{\mathcal O}'\subset {\mathcal O}$ and 
$|\pi_\bq(\theta)|>1$ for $\theta\in {\mathcal O}'$. For each 
$\varepsilon>0$ consider the function 
$$ 
\alpha_\varepsilon(\theta)=v(\theta)+\log\frac{2}{1+|\pi_\bq(\theta)|^2}+C_0+ 
(2+\varepsilon)\log |\pi_\bq(\theta)|, \quad \theta\in 
\overline{\mathcal O}'\backslash\{\bq\}, 
$$ 
where $C_0$ is a constant such that 
$v(\theta)+\log\frac{2}{1+|\pi_\bq(\theta)|^2}+C_0>0$ on $\p 
{\mathcal O}'$. One can check that 
$$ 
\Delta_{{\mathbb S}^2}\left(\log 
\frac{2}{1+|\pi_\bq(\theta)|^2}\right)=-1 \quad \mbox{and} \quad 
\Delta_{{\mathbb S}^2} (\log|\pi_\bq(\theta)|)=0 \quad \mbox{in } 
{\mathcal O}'\setminus\{\bq\}. 
$$ 
Moreover $\lim_{\theta\rightarrow \bq} |\pi_\bq(\theta)|=+\infty$. 
Noting that (\ref{3.32}) implies 
$$ 
\liminf_{\theta\in {\mathcal O}, \theta\rightarrow \bq} 
\frac{v(\theta)}{\log |\pi_\bq(\theta)|}\ge 0. 
$$ 
We therefore have 
$$ 
-\Delta_{{\mathbb S}^2}\alpha_\varepsilon>0 \,\mbox{ in } {\mathcal 
O}'\backslash \{\bq\}, \,\, \alpha_\varepsilon>0 \mbox{ on } \p 
{\mathcal O}'\, \, \mbox{ and } \,\, \lim_{\theta\in {\mathcal O}, 
\theta\rightarrow \bq} \alpha_\varepsilon=+\infty. 
$$ 
It then follows from the maximum principle that 
$\alpha_\varepsilon>0$ on $\overline{\mathcal O}'\setminus \{\bq\}$. 
Letting $\varepsilon\rightarrow 0$ we get 
$$ 
v(\theta)\ge -C_0-\log 
\frac{2|\pi_\bq(\theta)|^2}{1+|\pi_\bq(\theta)|^2}\ge -C_0-\log 2 
\quad \mbox{for } \theta\in {\mathcal O}'\setminus \{\bq\} 
$$ 
which gives the desired assertion. 
\end{proof}

\begin{Lemma}\label{L3.7} 
Let $\mathcal O$ be an open set in ${\mathbb S}^2$ and 
$\bq\in\overline{\mathcal O}$. If $v\in C^2({\mathcal 
O}\setminus\{\bq\})\cap C^0(\overline{\mathcal O}\backslash\{\bq\})$ 
is a non-negative function such that $-\Delta_{{\mathbb S}^2}v 
>0$ in ${\mathcal O}\setminus\{\bq\}$ and $v\ge c_0$ on 
$\p \mathcal{O} \backslash\{\bq\}$ for some constant $c_0>0$, then 
$v>c_1$ on $\overline{\mathcal O}\setminus\{\bq\}$ for some constant 
$c_1>0$. 
\end{Lemma} 
 
\begin{proof} First by the strong maximum principle we have $v>0$ 
on $\overline{O}\backslash\{\bq\}$. By shrinking $\mathcal {O}$ if 
necessary, we may assume $|\pi_\bq(\theta)|>1$ for $\theta\in 
{\mathcal O}$. For each $\varepsilon>0$ consider the function 
$$ 
\beta_\varepsilon(\theta):=v(\theta)-c_0+\varepsilon 
\log|\pi_\bq(\theta)|, \quad \theta\in \overline{\mathcal 
O}\backslash\{\bq\}. 
$$ 
Note that $\Delta_{{\mathbb S}^2}(\log|\pi_\bq(\theta)|)=0$ on 
$\mathcal{O}\backslash\{\bq\}$, $\lim_{\theta\rightarrow \bq} 
|\pi_\bq(\theta)|=+\infty$ and $v\ge c_0>0$ on $\p 
\mathcal{O}\backslash\{\bq\}$, we have 
$$ 
-\Delta_{{\mathbb S}^2} \beta_\varepsilon>0 \mbox{ in } {\mathcal 
O}\backslash\{\bq\},\,\, \, \beta_\varepsilon>0 \mbox{ on } \p 
{\mathcal O}\backslash\{\bq\}\,\, \mbox{ and } \,\,\lim_{\theta\in 
{\mathcal O}, \theta\rightarrow \bq} 
\beta_\varepsilon(\theta)=+\infty. 
$$ 
Therefore, by the maximum principle we have $\beta_\varepsilon>0$ 
on $\overline{\mathcal O}\backslash\{\bq\}$. Letting 
$\varepsilon\rightarrow 0$ gives the desired conclusion. 
\end{proof}

Now we are ready to give the proofs of Theorem \ref{T1.6} and Theorem \ref{T1.7}.

\begin{proof}[Proof of Theorem \ref{T1.6}] 
Since $f\ge 0$ and $K\ge 0$, we have 
$$ 
-\Delta_{{\mathbb S}^2} v+1>0 \quad \mbox{on } {\mathbb S}^2\backslash \{\bn, {\mathbf s}\}. 
$$ 
It then follows from Lemma \ref{L3.6} that 
$$ 
v\ge -C_0 \quad \mbox{on } {\mathbb S}^2\backslash\{\bn, {\mathbf s}\} 
$$ 
for some constant $C_0>0$. Lemma \ref{L3.5} then implies for each $\bp\in \p 
B_{\pi/2}({\mathbf s})$ there exists $0<\lambda_\bp<\pi/2$ such 
that 
$$ 
v_{\bp, \lambda}\le v \quad \mbox{on } \Sigma_{\bp, \lambda}\backslash\{\bn, {\mathbf s}\} 
\mbox{ for each } 0<\lambda<\lambda_\bp. 
$$ 
For each $\bp\in \p B_{\pi/2}({\mathbf 
s})$ we can define $\bar{\lambda}_\bp$ as before, then 
$\bar{\lambda}_\bp>0$. Using (\ref{1.16}), (\ref{3.31}), (K1), (f1), (Kf1), the symmetry 
properties of $K$ and $f$, and Lemma \ref{L3.7}, we may imitate the proof of Lemma \ref{L3.3} to 
conclude that 
$$ 
\bar{\lambda}_\bp\ge \pi/2 \quad \mbox{for each } \bp\in \p 
B_{\pi/2}({\mathbf s}). 
$$ 
 The desired assertion thus follows. 
\end{proof} 
 
\begin{proof}[Proof of Theorem \ref{T1.7}] 
Under the conditions in Theorem \ref{T1.7}, we can argue as in the 
proof of Theorem \ref{T1.6} to show that 
$$ 
v-v_{{\mathbf s},\pi/2}\ge 0 \quad \mbox{on } \Sigma_{{\mathbf s},\pi/2}\setminus\{\mathbf 
n\}. 
$$ 
Together with condition (K2) and (f2), we obtain that 
$$ 
-\Delta_{{\mathbb S}^2}(v-v_{{\mathbf s},\pi/2})\ge 
0,\quad\mbox{on}\quad \Sigma_{{\mathbf s},\pi/2}\setminus\{\mathbf 
n\}, 
$$ 
and, due to (Kf2), this inequality is strict somewhere in $\Sigma_{{\mathbf s},\pi/2}\setminus\{\mathbf 
n\}$. Hence we can follow the proof of Theorem \ref{T1.4} to show 
that 
$$ 
v-v_{{\bp},\pi/2}\ge 0,\quad\mbox{on}\quad 
\Sigma_{\bp,\pi/2}\setminus\{\mathbf n\} 
$$ 
for any $\bp \in \p B_{\pi/2}({\mathbf s})$. The proof 
is complete. 
\end{proof}

\section{\bf Appendix} 
\setcounter{equation}{0} 
 
In this section, we prove a Lemma from which Lemma \ref{L3.2} follows. For $n\ge 3$, let $\mathcal O$ be an open set in 
${\mathbb R}^n\setminus B_1(0)$. Consider 
 
\begin{equation}\label{A1} 
-\Delta 
v(y)+\sum_{i=1}^n\frac{b_i(y)}{|y|^3}v_i(y)+\frac{c(y)}{|y|^4}v(y)\ge 
0,\quad v>0,\quad\mbox{in}\quad {\mathcal O}, 
\end{equation}where $b_i, c\in L^{\infty}({\mathcal O})$. 
 
\begin{Lemma}\label{L4.1} 
For $n\ge 3$, let $\mathcal O$ be an open set in 
${\mathbb R}^n\setminus B_1(0)$. Assume $v\in C^2({\mathcal O})\cap C^0(\bar{\mathcal O})$ is a solution of (\ref{A1}). 
If $b_i, c\in L^{\infty}({\mathcal O})$ and there exists a constant $c_0>0$ such that 
$$ 
v(y)\ge \frac{c_0}{|y|^{n-2}}, \quad\mbox{on}\quad \p {\mathcal O} 
$$ and 
$$ 
\liminf_{|y|\to\infty}v(y)\ge 0. 
$$Then there exists a constant $c_1>0$ such that 
$$ 
v(y)\ge \frac{c_1}{|y|^{n-2}},\quad\mbox{in}\quad {\mathcal O}. 
$$ 
\end{Lemma} 
\begin{proof} The proof is based on the argument of  \cite[Lemma 2.1]{LZ}.  Let $\xi(y)=|y|^{2-n}+|y|^{1-n}$, then 
$$ 
-\Delta \xi(y)=-(n-1)|y|^{-n-1}. 
$$ 
By the condition on $b_i(x)$ and $c(x)$, there exists 
$R>1$ large enough, such that 
$$ 
-\Delta\xi(y)+\sum_{i=1}^n \frac{b_i(y)}{|y|^3}\xi_i(y)+\frac{c(y)}{|y|^4}\xi(y)\le 
0,\quad\mbox{in}\quad {\mathcal O}\setminus B_R(0). 
$$ 
Since $v>0$ in $\bar{\mathcal O}$, there exists some $\epsilon>0$ such 
that 
$$ 
v(y)\ge \epsilon\xi(y), \quad\mbox{for}\quad |y|=R,\quad y\in \bar{\mathcal O}, 
$$and 
$$ 
v(y)\ge \frac{c_0}{|y|^{2-n}}\ge \epsilon\xi(y),\quad\mbox{on}\quad 
\p{\mathcal O}. 
$$ 
We have 
$$ 
\left\{\begin{array}{l} 
        -\Delta(v-\epsilon\xi)+\sum\frac{b_i}{|y|^3}(v_i-\epsilon\xi_i)+\frac{c}{|y|^4}(v-\epsilon\xi)\ge 
        0, \quad {\mathcal O}\setminus B_R(0)\\ 
        \\ 
        v-\epsilon\xi\ge 0,\quad \mbox{on}\quad \p({\mathcal 
        O}\setminus B_R(0))\\ 
        \\ 
        \liminf_{|y|\to\infty}(v-\epsilon\xi)\ge 0. 
        \end{array}\right. 
$$ 
By the maximum principle, $v-\epsilon\xi\ge 0$ in ${\mathcal 
O}\setminus B_R(0)$. Thus for some constant $c_1>0$, 
$$ 
v\ge \frac{c_1}{|y|^{2-n}}, \quad\mbox{in}\quad {\mathcal O}. 
$$ 
\end{proof} 
 
We have the following equivalent lemma on a bounded open set of ${\mathbb R}^n$. 
For $n\ge 3$, let $\Omega$ be a bounded open set in ${\mathbb R}^n$, let $\bp\in\bar\Omega$. Consider 
\begin{equation}\label{A2} 
-\Delta v(y)+|y-\bp|\sum_{i=1}^n b_i(y)v_i(y)+c(y)v(y)\ge 0,\quad v>0,\quad\mbox{in}\quad \Omega, 
\end{equation}where $b_i,c\in L^{\infty}(\Omega)$. 
 
\begin{Lemma}\label{L4.2} 
For $n\ge 3$, let $\Omega$ be a bounded open set in ${\mathbb R}^n$, $\bp\in\bar\Omega$. Assume that $v\in C^2(\Omega\setminus\{\bp\})\cap C^0(\bar\Omega\setminus\{\bp\})$ is a solution of (\ref{A2}). 
If $b_i, c\in L^{\infty}(\Omega)$ and 
$$ 
\inf_{\p\Omega\setminus\{\bp\}}v>0. 
$$ 
Then $v>c_1$ in $\Omega\setminus\{\bp\}$ for some constant $c_1>0$. 
\end{Lemma} 
\begin{proof} 
Let 
$$ 
\tilde{\Omega}:=\left\{z\in {\mathbb R}^n|\bp+\frac{z-\bp}{|z-\bp|^2}\in\Omega\right\}, 
$$ 
and make a Kelvin transform with respect to $\bp$, 
$$ 
\tilde{v}(z)=\frac{1}{|z-\bp|^{n-2}}v\left(\bp+\frac{z-\bp}{|z-\bp|^2}\right). 
$$ 
Then $\tilde{v}\in  C^2(\tilde{\Omega})\cap C^0(\bar{\tilde{\Omega}})$. It is easy to check that $\tilde{v}$ satisfies 
$$ 
-\Delta\tilde{v}+\sum_{i=1}^n\frac{\tilde{b}_i(z)}{|z-\bp|^3}\tilde{v}_i(z)+\frac{\tilde{c}(z)}{|z-\bp|^4}\tilde{v}(z)\ge 0, \quad\tilde{v}>0,\quad\mbox{in}\quad \tilde{\Omega}, 
$$ 
where $\tilde{b}_i,\tilde{c}\in L^{\infty}(\tilde{\Omega})$ and 
$$ 
\tilde{v}(z)\ge \frac{c_0}{|z-\bp|^{n-2}},\quad\mbox{on}\quad\p 
\tilde{\Omega}. 
$$for some constant $c_0>0$. 
>From Lemma \ref{L4.1}, we get $\tilde{v}(z)\ge \frac{c_1}{|z-\bp|^{n-2}}$ 
in $\tilde{\Omega}$ for some constant $c_1>0$. The conclusion of Lemma follows easily. 
\end{proof} 
 
\begin{Remark}\label{R.A.1} 
It is easy to see from the proof that $\Delta$ can be replace by $a_{ij}\p_{ij}$ with $a_{ij}\in C^0(\bar\Omega)$, $(a_{ij})>0$ in $\bar{\Omega}$. In fact the same conclusion holds for more general operators. 
\end{Remark} 
 
\section*{\bf Acknowledgement} Q. Jin would like to thank MSRI for the financial 
support during the visit in the fall of 2006. Part of the work was completed while 
Y. Y. Li was visiting D\'{e}partement de Math\'{e}matique, Universit\'{e} Paris VI, 
he thanks the host T. Aubin and the institute for the kind invitation and hospitality.

\end{document}